\documentclass[12pt,a4paper]{amsart}

\usepackage{amssymb}
\usepackage{amscd}
\usepackage{hyperref}
\usepackage[numeric, abbrev, nobysame]{amsrefs}

\def\frak{\mathfrak}
\def\Bbb{\mathbb}
\def\Cal{\mathcal}

\let\phi\varphi

\newcommand{\x}{\times}
\renewcommand{\o}{\circ}

\newcommand{\al}{\alpha}
\newcommand{\be}{\beta}
\newcommand{\ga}{\gamma}

\newcommand{\ep}{\epsilon}

\newcommand{\la}{\lambda}
\newcommand{\om}{\omega}
\newcommand{\ph}{\phi}
\newcommand{\ps}{\psi}
\renewcommand{\th}{\theta}
\newcommand{\si}{\sigma}

\newcommand{\Ga}{\Gamma}
\newcommand{\La}{\Lambda}
\newcommand{\Ph}{\Phi}

\newcommand{\Om}{\Omega}

\newcommand{\Ups}{\Upsilon}

\newcommand{\im}{\operatorname{im}}

\newcommand{\id}{\operatorname{id}}

\newcommand{\Ad}{\operatorname{Ad}}
\newcommand{\Aut}{\operatorname{Aut}}
\newcommand{\gr}{\operatorname{gr}}
\newcommand{\pr}{\operatorname{pr}}
\newcommand{\Fl}{\operatorname{Fl}}

\newcommand{\tcg}{{\tilde{\Cal G}}}

\newcounter{theorem}
\newcounter{proposition}
\newcounter{lem}
\newcounter{corollary}
\newtheorem{thm}[theorem]{Theorem}
\newtheorem*{thm*}{Theorem \thesubsection}
\newtheorem{lemma}[lem]{Lemma}
\newtheorem{prop}[proposition]{Proposition}
\newtheorem{cor}[corollary]{Corollary}
\newtheorem*{lemma*}{Lemma \thesubsection}
\newtheorem*{prop*}{Proposition \thesubsection}
\newtheorem*{cor*}{Corollary \thesubsection}

\newcounter{def}
\newcounter{rem}
\theoremstyle{definition}
\newtheorem{definition}[def]{Definition}
\newtheorem*{definition*}{Definition \thesubsection}

\newtheorem*{example*}{Example \thesubsection}
\theoremstyle{remark}
\newtheorem{remark}[rem]{Remark}
\newtheorem*{remark*}{Remark \thesubsection}

\def\sideremark#1{\ifvmode\leavevmode\fi\vadjust{\vbox to0pt{\vss% the remark
 \hbox to 0pt{\hskip\hsize\hskip1em%                          will appear only
 \vbox{\hsize3cm\tiny\raggedright\pretolerance10000%          on the side
  \noindent #1\hfill}\hss}\vbox to8pt{\vfil}\vss}}}%
\begin{document}

\title{Parabolic conformally symplectic\\ 
structures II; parabolic contactification} 
\date{October 3, 2017}
\author{Andreas \v Cap and Tom\'a\v s Sala\v c}
\thanks{Support by projects P23244--N13 (both authors) and P27072--N25 
  (first author) of the Austrian Science fund (FWF) is gratefully
  acknowledged. We thank the referee for very helpful comments.}

\address{A.\v C.: Faculty of Mathematics\\
University of Vienna\\
Oskar--Morgenstern--Platz 1\\
1090 Wien\\
Austria}
\address{T.S.: Mathematical Institute\\ Charles University\\ Sokolovsk\'a
  83\\Praha\\Czech Republic}
\email{Andreas.Cap@univie.ac.at}
\email{salac@karlin.mff.cuni.cz}

\begin{abstract}
  Parabolic almost conformally symplectic structures were introduced
  in the first part of this series of articles as a class of geometric
  structures which have an underlying almost conformally symplectic
  structure.  If this underlying structure is conformally symplectic,
  then one obtains a PCS--structure. In the current article, we relate
  PCS--structures to parabolic contact structures. Starting from a
  parabolic contact structure with a transversal infinitesimal
  automorphism, we first construct a natural PCS--structure on any
  local leaf space of the corresponding foliation. Then we develop a
  parabolic version of contactification to show that any
  PCS--structure can be locally realized (uniquely up to isomorphism)
  in this way.

  In the second part of the paper, these results are extended to an
  analogous correspondence between contact projective structures and
  so--called conformally Fedosov structures.  The developments in this
  article provide the technical background for a construction of
  sequences and complexes of differential operators which are naturally
  associated to PCS--structures by pushing down BGG sequences on
  parabolic contact structures. This is the topic of the third part of
  this series of articles.
\end{abstract}

\subjclass[2010]{Primary: 53D05, 53D10, 53B15, 53C10, 53C15; 
Secondary: 53C29, 53C55}

\keywords{almost conformally symplectic structure, special symplectic
  connection, parabolic contact structure, conformally Fedosov
  structure, contactification}

\maketitle

\pagestyle{myheadings} \markboth{\v Cap and Sala\v c}{PCS--structures
  II}
 
\section{Introduction}\label{1}
This is the second part in a series of three articles devoted to the
study of a class of geometric structures and of differential complexes
naturally associated to these structures. The main motivation for this
series are the examples of differential complexes on complex
projective space, which were constructed and applied to questions of
integral geometry in \cite{E-G}. An attempt to put these complexes into
the context of geometric structures was made in the first version of
the preprint \cite{E-S}, which introduced the concept of a conformally
Fedosov structure. The aim of providing a framework analogous to
Bernstein--Gelfand--Gelfand resolutions (or BGG resolutions)
associated to parabolic geometries was realized (for conformally
Fedosov structures) in the second version of \cite{E-S}, which has
appeared recently.

Our article builds on \cite{PCS1}, where we introduced a family of
first order structures, which all have an underlying almost
conformally symplectic structure. There is one such structure for each
contact grading of a simple Lie algebra, which is not of type
$C_n$. These gradings are related to certain parabolic subalgebras and
to parabolic contact structures as discussed below, which motivates
the name \textit{parabolic almost conformally symplectic structures}
(or PACS--structures for short). The main result of \cite{PCS1} is
that each such structure determines a canonical connection on the
tangent bundle, which is characterized by a normalization condition on
its torsion. The torsion of this connection is a basic invariant of
the structure, which naturally splits into two parts. One of these
parts is the obstruction to the almost conformally symplectic
structure being conformally symplectic, and requiring this part of the
torsion to vanish, one arrives at the subclass of PCS--structures or
parabolic conformally symplectic structures.

In this second part, we only consider PCS--structures and we relate
them to another class of geometric structures, called parabolic
contact structures. There is one such structure for each contact
grading of a simple Lie algebra, and they all have an underlying
contact structure. The most prominent example of a parabolic contact
structure is provided by (partially integrable almost) CR--structures
of hypersurface type. Via the relation to contact gradings, each type
of PACS--structure determines a corresponding type of parabolic contact
structure in one higher dimension.

Now for parabolic contact structures, there is the concept of
transversal infinitesimal automorphisms, which in particular are
transversal infinitesimal contactomorphisms for the underlying contact
structure. Any transversal infinitesimal contactomorphism defines a
foliation with one--dimensional leaves. In \cite{Cap-Salac}, we have
shown that any local leaf--space for such a foliation naturally
inherits a conformally symplectic structure. Further we have shown
that locally any conformally symplectic structure can be realized in
this way and that these realizations are locally unique up to
contactomorphism.

In the current article, we extend all these results to the setting of
PCS--structures and parabolic contact structures, see Theorems
\ref{thm2.5}, \ref{thm2.7}, and \ref{thm2.8}. Moreover, we discuss the
relation between distinguished connections associated to parabolic
contact structures and the canonical connections of
PCS--structures. Together with Section 4.7 of \cite{PCS1}, this allows
us to complete the discussion of the relation of PCS--structures to
special symplectic connections in the sense of
\cite{Cahen-Schwachhoefer} and thus to exceptional holonomies, see
Theorem \ref{thm2.10}.

The contact gradings of Lie algebras of type $C_n$ do not give rise to
a PACS--structure. However, there is a type of parabolic contact
structures associated to these gradings, the so--called contact
projective structures. For the applications in \cite{E-G} it was only
necessary to construct the initial parts of certain differential
complexes. However, the constructions can be extended, exhibiting that
the results look similar to the BGG resolutions associated to locally
flat contact projective structures. Also, the tractor bundle
associated to a conformally Fedosov structure in \cite{E-S} looks
similar to the standard tractor bundle of a contact projective
manifold in one higher dimension. So all that suggests that
conformally Fedosov structures might be related to parabolic contact
structures. In Section \ref{3} of this article, we show that this is
indeed the case and extend the above mentioned results on structures
on quotients and contactifications from \cite{Cap-Salac} also to this
setting. This will allow us to include the $C_n$--types into a uniform
treatment in the last part of this series.

There is also a class of special symplectic connections associated to
Lie algebras of type $C_n$, the so--called symplectic connections of
Ricci type, see \cite{Cahen-Schwachhoefer} and \cite{BC}. We also get
a characterization of these connections as conformally Fedosov
structures with locally flat parabolic contactification, see Corollary
\ref{cor3.3}. Together with Theorem \ref{thm2.10} this recovers the fact
that all special symplectic connections can be obtained from local
quotients of the homogeneous models of parabolic contact structures,
which is the crucial result of \cite{Cahen-Schwachhoefer}.

A third important aim of this article is to give examples of global
leaf spaces of compact contact manifolds which inherit PCS--structures
respectively conformally Fedosov structures. From the point of view of
differential complexes, such global quotients are interesting, since
they lead to information on the cohomology of the resulting sequences,
which are important for the applications in \cite{E-G}. As one of
these examples, we show that the Hopf--fibration $S^{2n+1}\to\Bbb
CP^n$ can be considered as a global parabolic contactification in two
ways. On the one hand, this concerns the CR--structure on $S^{2n+1}$
coming from the embedding as a hypersurface in $\Bbb C^{n+1}$, which
induces the K\"ahler structure on $\Bbb CP^n$ as the underlying
PCS--structure. On the other hand, viewing $S^{2n+1}$ as the space of
real rays in $\Bbb C^{n+1}$ (viewed as a real symplectic vector
space), it inherits a contact projective structure. This gives rise to
the conformally Fedosov structure on $\Bbb CP^n$ studied in
\cite{E-S}.

The second example of such a global space of leaves we discuss is a
bit more exotic. In Theorem \ref{thm2.6}, we show that the Grassmannian
of complex planes in $\Bbb C^{n+1}$ can be realized as a global space
of leaves for a transversal infinitesimal automorphism of the space of
those quaternionic lines in $\Bbb H^{n+1}$, which are isotropic for a
quaternionic skew--Hermitian form. This induces on the Grassmannian
the PCS--structure of quaternionic type from Corollary 3.8 of
\cite{PCS1}, which is homogeneous under $SU(n+1)$ and can be viewed as
a counterpart of the quaternion--K\"ahler metric on the Grassmannian
of planes.

In the last part \cite{PCS3} of the series, we show how invariant
differential operators for parabolic contact structures induce natural
differential operators on PCS--quotients. There is a well developed
theory which in particular produces invariant differential complexes
on locally flat parabolic contact structures and in certain curved
situations. These can then be used to obtain invariant differential
complexes associated to special symplectic connections as well as
certain more general PCS--structures.

\section{Contactification of PCS--structures}\label{2}
We start this section by briefly recalling the concept of
PCS--structures as introduced in \cite{PCS1} and the class of
parabolic contact structures. Then we describe the relation between
the two kinds of structures, generalizing the analogous results from
\cite{Cap-Salac}, which relate conformally symplectic structures to
contact structures.

\subsection{Contact gradings and associated groups}\label{2.1}
A \textit{contact grading} on a simple Lie algebra $\frak g$ over
$\Bbb K=\Bbb R$ or $\Bbb C$ is a decomposition $\frak g=\frak
g_{-2}\oplus\frak g_{-1}\oplus\frak g_0\oplus\frak g_1\oplus\frak g_2$
which is compatible with the Lie bracket in the sense that $[\frak
g_i,\frak g_i]\subset\frak g_{i+j}$ (putting $\frak g_\ell=\{0\}$ if
$|\ell|>2$), such that $\dim(\frak g_{-2})=1$ and such that the bracket
$\frak g_{-1}\x\frak g_{-1}\to\frak g_{-2}$ is a non--degenerate
bilinear form. It turns out that then the analogous statements hold
for $\frak g_2$ and the bracket $\frak g_1\x\frak g_1\to\frak g_2$. 

Such gradings are well known from the theory of quaternionic symmetric
spaces as well as from the theory of parabolic contact
structures. They exist on any complex simple Lie algebra of rank at
least two and on most non--compact real simple Lie algebras (namely on
those which contain a highest root vector). If they exist, they are
always unique up to isomorphism. The complete list of contact gradings
can be found in Example 3.2.10 of \cite{book}.

Via $\frak g_0\subset\frak g$, which is a Lie subalgebra by the
grading property, contact gradings give rise to geometric structures
in two ways. On the one hand, consider the graded Lie algebra $\frak
g_-:=\frak g_{-2}\oplus\frak g_{-1}$ which, by definition of a contact
grading, is a Heisenberg algebra. By the grading property and the
Jacobi--identity, the adjoint action of any element $A\in\frak g_0$
can be restricted to $\frak g_-$ and defines a derivation on this Lie
algebra which preserves the grading. It is easy to see that the
resulting Lie algebra homomorphism $\frak g_0\to\frak{der}_{\gr}(\frak g_-)$
is always injective.

It is easy to describe $\frak{der}_{\gr}(\frak g_-)$
explicitly. Viewing the bracket as a linear map $\La^2\frak
g_{-1}\to\frak g_{-2}$, its kernel defines a codimension one--subspace
$\La^2_0\frak g_{-1}\subset\La^2\frak g_{-1}$. Now one easily shows
that a linear map $\ph:\frak g_{-1}\to\frak g_{-1}$ can be extended to
a derivation of $\frak g_-$ if and only if the induced endomorphism
of $\La^2\frak g_{-1}$ preserves the hyperplane $\La^2_0\frak g_{-1}$
and thus factors to an endomorphism of the quotient which is
isomorphic to $\frak g_{-2}$. The bracket defines a non--degenerate
bilinear form on $\frak g_{-1}$ up to scale and thus a line in
$\La^2(\frak g_{-1})^*$, which is exactly the annihilator of
$\La^2_0\frak g_{-1}$. The stabilizer of this line in $L(\frak
g_{-1},\frak g_{-1})$ will be denoted by $\frak{csp}(\frak g_{-1})$
and called the \textit{conformally symplectic Lie algebra} of $\frak
g_{-1}$. The above discussion implies that the adjoint action of
$\frak g_0$ on $\frak g_{-1}$ defines an injective homomorphism $\frak
g_0\hookrightarrow\frak{csp}(\frak g_{-1})$.

The case of algebras of type $C_n$ is special with respect to contact
gradings. In this case $\frak g$ itself is a symplectic Lie algebra
and it turns out that the adjoint action actually gives rise to an
isomorphism $\frak g_0\cong \frak{csp}(\frak
g_{-1})\cong\frak{der}_{gr}(\frak g_-)$, so one just obtains the full
conformally symplectic algebra in this case. Therefore, we will have
to exclude the algebras of type $C_n$ from the considerations in the
rest of Section \ref{2}. Conformally Fedosov structures and their
contactifications, which will be discussed in Section \ref{3} below,
can be thought of as a $C_n$--analog of the situation discussed here.

To proceed towards geometric structures determined by a contact
grading, we also need a choice of group associated to the Lie algebra
$\frak g_0$. While the theory of PACS--structures can be developed for
more general Lie groups with Lie algebra $\frak g_0$, we will restrict
the choice of group here to get a complete correspondence to parabolic
contact structures. Namely, let us start with a Lie group $G$ with Lie
algebra $\frak g$. Then it is well known that the Lie subalgebra
$\frak p:=\frak g_0\oplus\frak g_1\oplus\frak g_2\subset\frak g$ is a
parabolic subalgebra, so in particular, the normalizer of $\frak p$ in
$G$ is a Lie subgroup with Lie algebra $\frak p$. Let $P\subset G$ be
a \textit{parabolic subgroup} of $G$ corresponding to $\frak p$,
i.e.~a subgroup lying between this normalizer and its connected
component of the identity. Then $P$ has Lie algebra $\frak p$ and it
is well known that the adjoint action of any element of $P$ maps each
$\frak g_i$ to $\frak g_i\oplus\dots\oplus\frak g_2$, so it preserves
the filtration of $\frak g$ induced by the grading.

Finally, one defines a Lie subgroup $G_0\subset P$ as consisting of
those elements whose adjoint action preserves the grading on $\frak
g$. Let us denote by $CSp(\frak g_{-1})\subset GL(\frak g_{-1})$ the
subgroup of those linear isomorphisms for which the induced
automorphism of $\La^2(\frak g_{-1})^*$ preserves the line determined
by the bracket. Then this can be identified with the group
$\Aut_{gr}(\frak g_-)$ of automorphisms of the graded Lie algebra
$\frak g_-$ and the adjoint action defines an infinitesimally
injective homomorphism $G_0\to CSp(\frak g_{-1})$. 

\subsection{PCS--structures}\label{2.2} 
Having constructed the infinitesimally injective homomorphism $G_0\to
CSp(\frak g_{-1})\subset GL(\frak g_{-1})$ there is the natural notion
of a first order structure with structure group $G_0$ on smooth
manifolds of dimension $\dim(\frak g_{-1})$. These are the
PACS--structures as introduced in \cite{PCS1}. The simplest way to
describe such a structure is as a principal bundle $p:\Cal G_0\to M$
with structure group $G_0$, which is endowed with a strictly
horizontal, $G_0$--equivariant one--form $\th\in\Om^1(\Cal G_0,\frak
g_{-1})$. Here equivariancy means that for the principal right action
$r^g$ of $g\in G_0$, one has $(r^g)^*\th=\Ad(g^{-1})\o\th$ while
strict horizontality means that in each point $u\in\Cal G_0$, the
kernel of $\th(u):T_u\Cal G_0\to\frak g_{-1}$ coincides with the
vertical subbundle of $\Cal G_0\to M$. In particular, $\th(u)$ descends
to a linear isomorphism $T_{p(u)}M\to\frak g_{-1}$ which identifies
$\Cal G_0$ as a reduction of the linear frame bundle of $TM$. 

Given $(\Cal G_0\to M,\th)$ as above, a representation of $G_0$ on a
vector space $V$ gives rise to the natural vector bundle $\Cal
G_0\x_{G_0}V\to M$. Via $\th$, the bundle $\Cal G_0\x_{G_0}\frak
g_{-1}$ is identified with $TM$. Further, the $G_0$--invariant line in
$\La^2(\frak g_{-1})^*$ gives rise to a line subbundle
$\ell\subset\La^2T^*M$ such that each non--zero element of $\ell$ is
non--degenerate as a bilinear form on the corresponding tangent
space. This is the almost conformally symplectic structure underlying
a PACS--structure. While the theory of PACS--structures is developed
in this general setting in \cite{PCS1}, we will impose a restriction
on such structures throughout this article. Namely, we will only deal
with PCS--structures, i.e.~assume that the underlying structure is
conformally symplectic. This means that the line subbundle $\ell$
admits local sections which are closed as two--forms on $M$. 

A central result of \cite{PCS1} is that any PACS--structure on $M$
determines a canonical linear connection on $TM$, whose torsion
satisfies a certain normalization condition. We formulate this result
only for PCS--structures, where the normalization condition on the
torsion is simpler. To formulate this condition, we need the fact that
there a natural $G_0$--invariant subspace
$\ker(\square)\subset\La^2(\frak g_{-1})^*\otimes\frak g_{-1}$, called
the \textit{harmonic subspace}. Here $\square:\La^2(\frak
g_-)^*\otimes\frak g_-\to \La^2(\frak g_-)^*\otimes\frak g_-$ is the
Kostant--Laplacian, see Section 4.3 of \cite{PCS1} for the definition
and Sections 4.5 and 4.6 of this reference as well as Section 4.2 of
\cite{book} for an explicit description of this subspace in several
cases. Via associated bundles, this gives rise to a subbundle
$\ker(\square)\subset\La^2T^*M\otimes TM$, whose elements are called
\textit{algebraically harmonic}. Specialized to the PCS--case,
Corollary 4.3 of \cite{PCS1} shows
\begin{thm}\label{thm2.2}
  Let $(p:\Cal G_0\to M,\th)$ be a PCS--structure (associated to a
  contact grading on a simple Lie algebra which is not of type $C_n$)
  on a smooth manifold $M$. Then there is a unique connection
  compatible with this structure, such that the induced linear
  connection on $TM$ has algebraically harmonic torsion.
\end{thm}

\subsection{Parabolic contact structures}\label{2.3} 
There is a second way to obtain a geometric structure from a contact
grading of a real simple Lie algebra $\frak g$. Consider a smooth
manifold $M^\#$ of dimension $\dim(\frak g_-)=\dim(\frak g_{-1})+1$
and suppose that $H\subset TM^\#$ is a smooth distribution of corank
one. Denoting by $Q$ the quotient bundle $TM^\#/H$ (which is a real
line bundle), the Lie bracket of vector fields induces a skew
symmetric bilinear bundle map $\Cal L:H\x H\to Q$ called the
Levi--bracket. Now $H$ is a contact structure if and only if this
Levi--bracket is non--degenerate in each point $x\in
M^\#$. Equivalently, for each point $x$, the associated graded
$H_x\oplus Q_x$ of the tangent space $T_xM^\#$ endowed with $\Cal L_x$
as a Lie bracket is isomorphic to $\frak g_-$. Supposing that
$H\subset TM^\#$ is a contact structure, it is then easy to obtain an
adapted frame bundle for the associated graded vector bundle
$\gr(TM^\#)$ to the filtered vector bundle $H\subset TM^\#$ with
structure group $\Aut_{gr}(\frak g_-)\cong CSp(\frak g_{-1})$.

If $\frak g$ is not of type $C_n$, then a parabolic contact structure
of the corresponding type can be simply described as a reduction of
structure group of this adapted frame bundle corresponding to the
infinitesimally injective homomorphism $G_0\to \Aut_{gr}(\frak g_-)$
defined by the adjoint action. Similarly to the case of standard first
order structures, such a reduction can be described as an abstract
principal $G_0$--bundle $p_0^\#:\Cal G_0^\#\to M^\#$ endowed with an analog
of a soldering form as follows. 

Define $T^{-1}\Cal G_0^\#\subset T\Cal G_0^\#$ as the preimage of
$H\subset TM^\#$ under the natural projection. Of course, this
contains the vertical subbundle $V\Cal G_0^\#=\ker(Tp_0^\#)$. Then the
soldering form consists of two components, namely
$\th_{-2}^\#\in\Om^1(\Cal G_0^\#,\frak g_{-2})$ and
$\th_{-1}^\#\in\Ga(L(T^{-1}\Cal G_0^\#,\frak g_{-1}))$. So while
$\th_{-2}^\#$ is an ordinary differential form, $\th_{-1}^\#$ is only
defined on the subbundle $T^{-1}\Cal G_0^\#$. Both components are
required to be $G_0$--equivariant (with respect to the adjoint
action), and they should be strictly horizontal in the sense that for
each $u\in\Cal G_0^\#$ one has $\ker(\th_{-2}^\#(u))=T^{-1}_u\Cal
G_0^\#$ and $\ker(\th_{-1}^\#(u))=V_u\Cal G_0^\#$.  This implies that
$\th_{-2}^\#(u)$ descends to a linear isomorphism
$Q_{p_0^\#(u)}\to\frak g_{-2}$ while $\th_{-1}^\#(u)$ descends to a
linear isomorphism $H_{p_0^\#(u)}\to\frak g_{-1}$. These isomorphisms
together define an isomorphism of Lie algebras, i.e.~they intertwine
between the Levi--bracket and the Lie bracket on $\frak g_-$. 

Parabolic contact structures admit an equivalent uniform description,
which also works in the $C_n$--case. Namely, given a smooth manifold
$M^\#$ as above, one considers a principal bundle $p^\#:\Cal G^\#\to
M^\#$ with structure group $P$, which is endowed with a Cartan
connection $\om\in\Om^1(\Cal G^\#,\frak g)$. This by definition means
that $\om$ is $P$--equivariant and reproduces the generators of
fundamental vector fields, and that for each $u\in\Cal G^\#$ the map
$\om(u):T_u\Cal G^\#\to\frak g$ is a linear isomorphism. Moreover, one
has to assume the conditions of regularity and normality on the
curvature of the Cartan connection $\om$, see Section 3.1 of
\cite{book} for the precise definitions.

The two descriptions are related as follows.  Given the Cartan
geometry $(\Cal G^\#,\om)$ and a point $u\in\Cal G^\#$ with
$x=p^\#(u)\in M^\#$, $\om(u)$ descends to a linear isomorphism
$T_xM^\#\to \frak g/\frak p$. The subset $\frak g_{-1}\oplus\frak
p\subset \frak g$ is $P$--invariant so its preimage gives rise to a well
defined subspace $H_x\subset T_xM^\#$. These subspaces fit together to
a corank one subbundle $H\subset TM^\#$ which defines a contact structure
by regularity of $\om$. Moreover, the exponential mapping restricts to
a diffeomorphism from $\frak p_+:=\frak g_1\oplus\frak g_2$ onto a
closed normal subgroup $P_+\subset P$, such that $G_0\cong P/P_+$ via
the restriction of the canonical projection. The group $P_+$ acts
freely on $\Cal G^\#$ via the principal right action, and the quotient
$\Cal G_0^\#:=\Cal G^\#/P_+$ naturally is a principal bundle over
$M^\#$ with structure group $P/P_+=G_0$. Finally, one shows that the
component of $\om$ in $\frak g_{-2}$ descends to
$\th_{-2}^\#\in\Om^1(\Cal G_0^\#,\frak g_{-2})$ while an appropriate
restriction of the $\frak g_{-1}$--component of $\om$ descends to
$\th_{-1}^\#\in\Ga(L(T^{-1}\Cal G_0^\#,\frak g_{-1}))$. 

In the $C_n$--case, the underlying structure constructed in this way
is just the adapted frame bundle of a contact structure on $M^\#$ and
does not contain any additional information. In all other cases,
however, the Cartan geometry $(\Cal G^\#,\om)$ can be reconstructed
from the underlying structure by a prolongation procedure as described
in Section 3.1 of \cite{book}. One forms the obvious extension $\Cal
G_0^\#\x_{G_0}P$ of structure group and then shows that there is a
normal Cartan connection on this principal bundle inducing the given
soldering form. Further, one proves that any isomorphism of underlying
structures lifts to an isomorphism of normal Cartan geometries, which
establishes an equivalence (in the categorical sense) between the two
pictures.

\subsection{PCS--quotients and parabolic contactifications}\label{2.4}  

Recall that for a contact manifold $(M^\#,H)$ an \textit{infinitesimal
  contactomorphism} is a vector field $\xi\in\frak X(M^\#)$ such that
for all $\eta\in\Ga(H)$ we get $[\xi,\eta]\in\Ga(H)$. An infinitesimal
contactomorphism is called \textit{transversal} if $\xi(x)\notin H_x$
for all $x\in M^\#$. This in particular implies that $\xi$ is nowhere
vanishing and thus defines a foliation of $M^\#$ with one--dimensional
leaves. By a \textit{quotient} of $M^\#$ by $\xi$ one then means a
global space of leaves, i.e.~a surjective submersion $q:M^\#\to M$
with connected fibers such that for each $x\in M^\#$ the vertical
subspace $\ker(T_xq)\subset T_xM^\#$ coincides with the line spanned
by $\xi(x)$.

In this situation, Proposition 2.2 of \cite{Cap-Salac} shows that the
contact structure on $M^\#$ induces a conformally symplectic structure
$\ell\subset\La^2 T^*M$. Denoting by $\al$ the unique contact form on
$M^\#$ such that $\al(\xi)=1$, there even is a global symplectic form
$\om$ on $M$, which is a section of $\ell$ and such that
$q^*\om=d\al$. Conversely, given a manifold $M$ endowed with a
conformally symplectic structure $\ell$, one can try to (locally)
realize it as a quotient of a contact manifold $(M^\#,H)$, which then
is called a (local) \textit{contactification} of $(M,\ell)$. Existence
and uniqueness of such contactifications is established in
\cite{Cap-Salac}.

We want to develop refinements of these concepts which relate parabolic
contact structures to PCS--structures. Let us start with the regular
normal parabolic geometry $(p:\Cal G^\#\to M^\#,\om)$ of type $(G,P)$
corresponding to a parabolic contact structure with contact
distribution $H\subset TM^\#$. An \textit{infinitesimal symmetry} of
the geometry then is a vector field $\tilde\xi\in\frak X(\Cal G^\#)$
which is $P$--invariant and satisfies $L_{\tilde\xi}\om=0$, where $L$
denotes the Lie derivative. A $P$--invariant vector field is
automatically projectable, so there is a corresponding vector field
$\xi\in\frak X(M^\#)$ as well as an intermediate $G_0$--invariant
vector field $\xi_0\in\frak X(\Cal G_0^\#)$. It is easy to verify
directly that $\xi$ is an infinitesimal contactomorphism. As before,
an infinitesimal symmetry is called \textit{transverse} if
$\xi(x)\notin H_x$ for all $x\in M^\#$.

In the case that $\frak g$ is not of type $C_n$, infinitesimal
symmetries can be equivalently described in terms of the underlying
structure $(p_0^\#:\Cal G_0^\#\to M^\#,\th^\#)$. For a
$G_0$--invariant vector field $\xi_0\in\frak X(\Cal G_0^\#)$ to be an
infinitesimal symmetry, one first has to require that
$[\xi_0,\eta]\in\Ga(T^{-1}\Cal G_0^\#)$ for all $\eta\in\Ga(T^{-1}\Cal
G_0^\#)$ or equivalently that the projection $\xi\in\frak
X(M^\#)$ of $\xi_0$ is an infinitesimal contactomorphism. If this is
the case, there is a well defined Lie derivative
$L_{\xi_0}\th^\#_{-1}$ of the partially defined differential form
$\th_{-1}^\#$, and $\xi_0$ is an infinitesimal symmetry if
$L_{\xi_0}\th^\#_i=0$ for $i=-1,-2$. 

\begin{definition}\label{def2.4}
  Consider a parabolic contact structure of type $(G,P)$ with
  underlying structure $(p_0^\#:\Cal G_0^\#\to M^\#,\th^\#)$ and a
  transversal infinitesimal symmetry $\xi_0\in\frak X(\Cal G_0^\#)$ of
  this structure. Then a \textit{PCS--quotient} of the parabolic
  contact structure is a PCS--structure $(p:\Cal G_0\to M,\th)$ with
  structure group $G_0$ together with a morphism $q_0:\Cal
  G_0^\#\to\Cal G_0$ of principal bundles such that
  \begin{itemize}
  \item $q_0$ is a surjective submersion with connected fibers. 
  \item For each $u^\#\in\Cal G_0^\#$ the kernel of $T_{u^\#}q_0$ is
    spanned by $\xi_0(u^\#)$.
  \item The restriction of $q_0^*\th\in\Om^1(\Cal G_0^\#,\frak
    g_{-1})$ to elements of $T^{-1}\Cal G_0^\#$ coincides with
    $\th^\#_{-1}$. 
  \end{itemize}
In this situation, $(p_0^\#:\Cal G_0^\#\to M^\#,\th^\#,\xi_0)$ is
also referred to as a \textit{parabolic contactification} of the
PCS--structure $(p:\Cal G_0\to M,\th)$. 
\end{definition}

To study these concepts, we need some preliminary
observations. Consider the underlying structure $(p_0^\#:\Cal
G_0^\#\to M^\#,\th^\#)$ of a parabolic contact structure of type
$(G,P)$, a transverse symmetry $\xi_0\in\frak X(\Cal G_0^\#)$ of this
structure and the projected vector field $\xi\in\frak X(M^\#)$. Then
both $\xi_0$ and $\xi$ are nowhere vanishing and hence define
foliations with one--dimensional leaves.

\begin{lemma}\label{lem2.4}
  In this situation, we have:

  (1) If $N\subset\Cal G_0^\#$ is a leaf of the foliation defined by
  $\xi_0$, then $p_0^\#(N)\subset M^\#$ is a leaf of the foliation
  defined by $\xi$ and the restriction of $p_0^\#$ to $N$ is a
  covering map $N\to p_0^\#(N)$.

  (2) Suppose that there is a PCS--structure $(p:\Cal G_0\to M,\th)$
  and a PCS--quotient $q_0:\Cal G_0^\#\to\Cal G_0$ by $\xi_0$. Then
  the base map $q:M^\#\to M$ of $q_0$ is a quotient by the transverse
  infinitesimal contactomorphism $\xi$. Moreover in this case the
  coverings from (1) actually have to be diffeomorphisms.
\end{lemma}
\begin{proof}
  (1) Since $\xi_0$ projects to $\xi$ and is transversal to the
  vertical subbundle of $p_0^\#$, it is clear that $p_0^\#$ maps small
  open subsets of $N$ to integral submanifolds of the distribution
  spanned by $\xi$. Hence $p_0^\#(N)$ is a connected immersed integral
  manifold for the foliation defined by $\xi$ and thus contained in
  some leaf $\underline{N}$ of this foliation. Again by construction,
  $p_0^\#|_N:N\to\underline{N}$ has bijective tangent maps in all
  points, so it is a local diffeomorphism and hence
  $p_0^\#(N)\subset\underline{N}$ is open. 

  Next, we observe that $\xi_0$ is $G_0$--equivariant. Hence if the
  flow of $\xi_0$ through a point $u^\#$ is defined on some time
  interval around zero, the same is true for $u^\#\cdot g$ for each
  $g\in G_0$. Now each such local flow line is either contained in $N$
  or disjoint from it. In particular, for $u^\#\in N$, such a local
  flow line projects to an open neighborhood of $p^\#_0(u^\#)$ in
  $\underline{N}$, whose pre--image in $N$ is the disjoint union of
  those flow lines, which are contained in $N$. By definition, these
  are open in the manifold topology of $N$.

  On the other hand, if $x^\#\in\underline{N}$ lies in the closure of
  $p^\#_0(N)$, then choose $u^\#\in\Cal G_0^\#$ lying over $x^\#$. For
  a local integral curve $c$ for $\xi_0$ through $u^\#$, $p_0^\#\o c$
  fills an open subset of $\underline{N}$ and hence intersects
  $p_0^\#(N)$. Shifting $c$ by the principal right action of an
  appropriate element of $G_0$, we may assume that $c$ intersects $N$,
  and hence is contained in $N$. Thus $x^\#\in p_0^\#(N)$, so
  $p_0^\#(N)\subset\underline{N}$ is closed and hence
  $p_0^\#(N)=\underline{N}$. Since we have constructed trivializing
  neighborhoods for $p_0^\#:N\to\underline{N}$ already, this completes
  the proof of (1).

\medskip

(2) By assumption, we have $p\o q_0=q\o p_0^\#$. Since both $p$ and
$q_0$ are surjective submersions, we see that $q$ is a surjective
submersion. Since $\xi_0$ projects to $\xi$, it is also clear that the
values of $\xi$ lie in $\ker(Tq)$ and since each of these kernels has
to be one--dimensional, it is spanned by the value of $\xi$. Now
consider a fiber $N$ of $q_0$. Then $p_0^\#(N)$ lies in one fiber of
$q$ and using $G_0$--equivariancy of $q_0$, one easily verifies that
$p_0^\#(N)$ is all of this fiber.  Hence the fibers of $q$ are
connected, so the first statement is proved. But now the fibers of
$q_0$ and $q$ are closed connected integral submanifolds for the
distributions defined by $\xi_0$ and $\xi$, and hence have to coincide
with the leaves. But since $q_0$ is a morphism of principal bundles,
its restriction to each fiber of $p_0^\#$ is injective, so $p_0^\#$ is
injective on fibers of $q_0$, and the last claim follows.
\end{proof}

\subsection{Existence of PCS--quotients}\label{2.5} 
Part (2) of Lemma \ref{lem2.4} tells us how to naturally phrase the
question of existence of PCS--quotients. We assume that we start with
a parabolic contact structure on $M^\#$ and the underlying vector
field $\xi\in\frak X(M^\#)$ of an infinitesimal symmetry of this
structure. Then the natural question to ask is when a quotient
$q:M^\#\to M$ of $M^\#$ by $\xi$ can be made into a PCS--quotient. We
can now show that the necessary condition from part (2) of Lemma
\ref{lem2.4} is also sufficient.

\begin{thm}\label{thm2.5}
  Consider a parabolic contact structure of type $(G,P)$ with $G$ not
  of type $C_n$ on a smooth manifold $M^\#$ with underlying
  $G_0$--bundle $p_0^\#:\Cal G_0^\#\to M^\#$. Let $\xi_0\in\frak
  X(\Cal G^\#_0)$ be a transverse infinitesimal symmetry of this
  structure and let $\xi\in\frak X(M^\#)$ be the underlying
  infinitesimal contactomorphism. Suppose that for each leaf
  $N\subset\Cal G_0^\#$ of the foliation defined by $\xi_0$, the
  restriction $p_0^\#|_N:N\to p_0^\#(N)$ is a diffeomorphism. 

  Then for any quotient $q:M^\#\to M$ by $\xi$, there is a canonical
  PCS--structure $\Cal G_0\to M$ on $M$ such that $q$ becomes a
  PCS--quotient.
\end{thm}
\begin{proof}
  Let $\sim$ be the equivalence relation on $\Cal G_0^\#$ defined by
  the foliation induced by $\xi_0$, i.e.~two points are equivalent if
  they lie in the same leaf. Define $\Cal G_0:=\Cal G_0^\#/\sim$, the
  set of equivalence classes, and let $q_0:\Cal G_0^\#\to\Cal G_0$ be
  the canonical map. From the proof of Lemma \ref{lem2.4}, we know that
  $p_0^\#$ maps leaves to leaves, and the leaves in $M^\#$ are the fibers
  of $q$, so there is a set map $p_0:\Cal G_0\to M$ such that $p_0\o
  q_0=q\o p_0^\#$. 

  We claim that $p_0:\Cal G_0\to M$ is a $G_0$--principal bundle and
  $q_0$ is a morphism of principal bundles. To prove this, observe
  first that $q\o p_0^\#$ is a surjective submersion. Hence for each
  $x\in M$, there is a neighborhood $U$ of $x$ in $M$ and a smooth map
  $\tau:U\to\Cal G_0^\#$ such that $q\o p_0^\#\o\tau=\id_U$. Using
  this, we define a map $\ps:U\x G_0\to\Cal G_0$ by
  $\ps(y,g):=q_0(\tau(y)\cdot g)$. This evidently satisfies
  $p_0\o\ps=\pr_1$, so it has values in $p_0^{-1}(U)$. If
  $\ps(y,g)=\ps(\tilde y,\tilde g)$, then applying $p_0$ we get
  $y=\tilde y$. Moreover, $\tau(y)\cdot g$ and $\tau(y)\cdot\tilde g$
  lie in the same leaf in $\Cal G_0^\#$. Since we have assumed that
  $p_0^\#$ restricts to an injection on each leaf, we conclude that
  $\tau(y)\cdot g=\tau(y)\cdot\tilde g$ and thus $g=\tilde g$. On the
  other hand, for $y\in U$, a point $z$ in $p_0^{-1}(U)$ corresponds
  to a leaf $N\subset\Cal G_0^\#$ such that
  $p_0^\#(N)=q^{-1}(\{y\})$. Hence there is a point $z^\#\in N$ such
  that $p_0^\#(z^\#)=p_0^\#(\tau(y))\in q^{-1}(\{y\})$. But this
  implies $z^\#=\tau(y)\cdot g$ for some $g\in G_0$ and hence
  $z=\ps(y,g)$, so $\ps:U\x G_0\to p_0^{-1}(U)$ is bijective.

  Now suppose that for some open subset $U\subset M$ we have two
  sections $\tau$ and $\hat\tau$ as above, and let us consider the
  ``chart change'' $\ph:U\x G_0\to G_0$, i.e.~the map characterized by
  $q_0(\tau(x)\cdot g)=q_0(\hat\tau(x)\cdot\ph(x,g))$ (which evidently
  exists). It suffices to do this locally around a fixed point $x$. By
  construction $p_0^\#(\tau(x))$ and $p_0^\#(\hat\tau(x))$ both lie in
  the leaf $q^{-1}(x)$, so there is a time $t_0\in\Bbb R$ such that
  $p_0^\#(\tau(x))=\Fl^{\xi}_{t_0}(p_0^\#(\hat\tau(x)))$. It is easy
  to see that the flow $\Fl^{\xi_0}_{t_0}$ is defined in the point
  $\hat\tau(x)$. Thus, we may shrink $U$ in such a way that for some
  $\ep>0$, $\Fl^{\xi_0}_{t_0+t}(\hat\tau(y))$ is defined for all $y\in
  U$ if $|t|<\ep$, and moreover $(y,t)\mapsto
  p_0^\#(\Fl^{\xi_0}_{t_0+t}(\hat\tau(y)))$ is a diffeomorphism from
  $U\x (-\ep,\ep)$ onto an open neighborhood $U^\#$ of
  $p_0^\#(\tau(x))$ in $M^\#$. Then there is an open neighborhood $V$
  of $x$ in $M$ such that $p_0^\#(\tau(V))\subset U^\#$, and we obtain
  a smooth function $\be:V\to (-\ep,\ep)$ such that
  $p_0^\#(\tau(y))=p_0^\#(\Fl^{\xi_0}_{t_0+\be(y)}\hat\tau(y))$.

  Now observe that
  $\tilde\tau(y):=\Fl^{\xi_0}_{t_0+\be(y)}(\hat\tau(y))$ is a smooth
  section of $q\o p_0^\#$ on $V$ which produces the same map $\ps$ as
  $\hat\tau$. But now $p_0^\#\o\tilde\tau=p_0^\#\o\tau$, so there is a
  smooth map $\ga:V\to G_0$ such that
  $\tilde\tau(y)=\tau(y)\cdot\ga(y)$, and hence
  $q_0(\tilde\tau(y)\cdot g)=q_0(\tau(y)\cdot (\ga(y)g))$. This shows
  that we can endow $\Cal G_0$ with a topology by requiring that the
  maps $\ps$ are homeomorphisms, and then there inverses define
  principal bundle charts, so the claim is proved.

  By construction, the quotient $q:M^\#\to M$ has the property that
  for each $x\in M^\#$, the tangent map $T_xq$ restricts to a linear
  isomorphism $H_x\to T_{q(x)}M$. Hence for each $u\in\Cal G_0^*$ the
  tangent map $T_uq_0$ restricts to a linear isomorphism $T^{-1}_u\Cal
  G_0^\#\to T_{q_0(u)}\Cal G_0$, which in addition respects the
  vertical subbundles. Composing $\th^\#_{-1}(u)$ with the inverse of
  this isomorphism, we obtain a surjective linear map $T_{q_0(u)}\Cal
  G_0\to \frak g_{-1}$ whose kernel is the vertical subspace
  $V_{q_0(u)}\Cal G_0$. Since $\th_{-1}^\#$ is preserved by the flow
  of $\xi_0$, it follows that locally around $u$ points in the same
  fiber lead to the same map on $T_{q_0(u)}\Cal G_0$. Since the fibers
  of $q_0$ are connected, we obtain a well defined map
  $\th(q_0(u)):T_{q_0(u)}\Cal G_0\to\frak g_{-1}$ with kernel the
  vertical subspace.
  
  Now $\ker(Tq_0)\subset T\Cal G_0^\#$ is a smooth subbundle which is
  complementary to $T^{-1}\Cal G_0^\#$, so we can use this to define a
  projection $\Pi:T\Cal G_0^\#\to T^{-1}\Cal G_0^\#$. Given
  $\xi\in\frak X(\Cal G_0)$, we can use a local smooth section $\si$
  of $q_0$ to write $\th(\xi)$ as $\th_{-1}^\#(\Pi(T\si\o\xi))$,
  which shows that $\th$ is smooth. Finally $\th$ is strictly
  horizontal by construction and $q_0^*\th$ restricts to $\th^\#_{-1}$
  on $T^{-1}\Cal G_0^\#$. Since equivariancy of $\th$ follows easily
  from equivariancy of $\th_{-1}^\#$, this completes the proof.
\end{proof}

\subsection{Examples}\label{2.6}
It is an easy consequence of the Frobenius theorem that given a
contact manifold $M^\#$ and a transversal infinitesimal
contactomorphism $\xi\in\frak X(M^\#)$, there locally exist quotients
by $\xi$ for which the fibers are intervals. Hence Theorem \ref{thm2.5}
shows that locally any transversal infinitesimal automorphism of a
parabolic contact structure admits PCS--quotients. Likewise, Theorem
\ref{thm2.5} shows that any global quotient $q:M^\#\to M$ can be made
into a PCS--quotient in case that all fibers of $q$ are simply
connected. We now discuss examples which show that in the case of
fibers which are circles, the situation is more subtle. Nonetheless,
we obtain several interesting examples of global contactifications in
which all fibers are circles.

These examples start from the homogeneous models of parabolic contact
structures. So we start with a simple Lie group $G$ and a parabolic
subgroup $P\subset G$ corresponding to a contact grading of the Lie
algebra $\frak g$ of $G$, and consider the homogeneous space $G/P$
with the parabolic contact structure coming from the Maurer--Cartan
form on $G$. The automorphism group of this geometry is $G$, so
looking for quotients, it is natural to consider the actions of
$1$--parameter subgroups of $G$ on $G/P$. We are particularly
interested in finding such actions for which the infinitesimal
generator is transversal everywhere, since these may lead to compact
quotients.

The first example concerns complex projective space $\Bbb CP^n$. The
Fubini--Study metric on $\Bbb CP^n$ is a K\"ahler metric, thus
defining a PCS--structure by Proposition 3.3 of \cite{PCS1}.

\begin{prop}\label{prop2.6}
  Consider $S^{2n+1}$ as the unit sphere in $\Bbb C^{n+1}$ and let
  $q:S^{2n+1}\to\Bbb CP^n$ be the Hopf fibration. Then endowing
  $S^{2n+1}$ with its usual CR--structure and $\Bbb CP^n$ with the
  PCS--structure of type $(PSU(n+1,1),P)$ defined the Fubini--Study
  metric, the map $q$ is a global parabolic contactification with
  circles as fibers.
\end{prop}
\begin{proof}
  Recall the description of $S^{2n+1}$ as the homogeneous model of
  strictly pseudoconvex CR--structures. Consider $V:=\Bbb C^{n+1}\x\Bbb
  C$ endowed with the standard Hermitian form of signature $(n+1,1)$,
  i.e.~the difference of the standard positive definite forms on the
  two factors.  Then any non--zero isotropic vector has to have
  non--zero first component. Mapping a point $z$ in the unit sphere
  $S^{2n+1}\subset\Bbb C^{n+1}$ to the line spanned by $(z,1)$
  identifies $S^{2n+1}$ with the space of isotropic complex lines in
  $V$. This gives rise to a transitive action of $G:=SU(n+1,1)$ on
  $S^{2n+1}$, thus identifying it with $G/P$, where $P$ is the
  stabilizer of an isotropic line. It is well known that the resulting
  diffeomorphisms of $S^{2n+1}$ are exactly those which preserve the
  strictly pseudoconvex CR--structure induced by the embedding
  $S^{2n+1}\hookrightarrow\Bbb C^{n+1}$.

Next, the obvious action of $U(1)$ on $S^{2n+1}$ by complex
multiplication can be realized by the action of a $1$--parameter
subgroup in $SU(n+1,1)$. Namely, for $t\in\Bbb R$, one considers
multiplication by $e^{it/(n+2)}$ in the first factor and
multiplication by $e^{-i(n+1)t/(n+2)}$ in the second factor. This
evidently defines a unitary map of determinant one and the point
$(e^{it/(n+2)}z,e^{-i(n+1)t/(n+2)})$ determines the same complex line
as $(e^{it}z,1)$. From this it is clear, that the infinitesimal
generator of this group is simply multiplication by $i$, which maps
any point into its real orthocomplement. However, $iz$ never lies in
the complex orthocomplement of $z$, which defines the CR
subspace. Hence we conclude that the infinitesimal automorphism
generating the one--parameter group is transversal everywhere, and
by definition the Hopf fibration $q$ is a global space of leaves for
this foliation. 

By construction, the $U(1)$--orbits of the action on $SU(n+1,1)$ are
$(n+2)$--fold coverings of the $U(1)$ orbits in $S^{2n+1}$, since the
projection sends $\lambda$ to multiplication by $\la^{n+2}$. The
underlying $G_0$--bundle in this case is simply $G/P_+\to G/P$, where
$P_+\subset P$ is the subgroup introduced in Section \ref{2.3}. It is
well known that its Lie algebra $\frak p_+$ consists of maps vanishing
on the distinguished complex line in $\Bbb C^{n+2}$, so any element of
$P_+$ acts as the identity on this line. This shows that $P_+$
intersects our one--parameter subgroup only in the identity, so the
projection $G\to G/P_+$ restricts to a diffeomorphism on each
orbit. Thus Lemma \ref{lem2.4} shows that $q:S^{2n+1}\to\Bbb CP^n$ cannot
be globally made into a PCS--quotient for the parabolic contact
structure of type $(G,P)$.

This is related to the fact that $G$ does not act effectively on
$G/P$, and correspondingly $G_0$ does not act effectively on $\Bbb
C^n$ (it is an $(n+2)$--fold covering of the conformal unitary group
$CU(n)$). This problem can be resolved by replacing $G$ by the
projective group $\underline{G}:=PSU(n+1,1)$ and $P$ by its image
$\underline{P}$ in $\underline{G}$. Then
$\underline{G}/\underline{P}=G/P$ and now the action becomes
effective. Moreover, $\underline{G}$ is the quotient of $G$ by its
center, which consists of the $(n+2)$nd roots of unity times the
identity map. This shows that the projection $G\to\underline{G}$
restricts to an $(n+2)$--fold covering on $U(1)$--orbits which exactly
identifies the different preimages of points in $S^{2n+1}$. Hence
Theorem \ref{thm2.5} implies that $q:S^{2n+1}\to\Bbb CP^n$ is globally a
PCS--quotient of geometries of type $(\underline{G},\underline{P})$,
for which $\underline{G}_0=CU(n)$.
\end{proof}

Our second example involves a type of parabolic contact structures
which has not been studied intensively in the literature, namely the
one associated to the groups $SO^*(2n)$. To formulate the necessary
background, consider a finite--dimensional right quaternionic vector
space $V$. Then a \textit{quaternionically skew Hermitian form} is a
map $\tau:V\x V\to\Bbb H$, which is bilinear over $\mathbb R$ and
satisfies $\tau(v,wq)=\tau(v,w)q$ and
$\tau(w,v)=-\overline{\tau(v,w)}$ for all $v,w\in V$ and $q\in\Bbb
H$. It is well known that any finite dimensional quaternionic vector
space admits a unique non--degenerate quaternionically skew Hermitian
form up to isomorphism. 

The group $SO^*(2n)$ is defined as the group of all quaternionically
linear automorphisms of $\Bbb H^n$ which preserve a non--degenerate
quaternionically skew Hermitian form $\tau$. It is easy to see that
$\tau$ can be recovered from its real part $\tau_{\Bbb R}$, so
preserving $\tau$ is equivalent to preserving $\tau_{\Bbb R}$. Of
course, $\tau_{\Bbb R}$ is just a skew symmetric bilinear form $\Bbb
H^n\x\Bbb H^n\to\Bbb R$ which satisfies $\tau_{\Bbb
  R}(va,wa)=\tau_{\Bbb R}(v,w)$ for any $a=i,j,k$ or equivalently for
any unit quaternion $a$.

For $n\geq 3$ consider the quaternionic projective space $\Bbb HP^n$
of one--dimensional quaternionic subspaces in $\Bbb H^{n+1}$. The
standard representation of $GL(n+1,\Bbb H)$ induces a transitive
action on $\Bbb HP^n$, which can be restricted to
$G:=SO^*(2n+2)$. However, this subgroup does not act transitively,
since the restriction of $\tau$ to a one--dimensional subspace can
either be zero or non--degenerate. It is easy to see that $G$ acts
transitively on the spaces of non--degenerate lines, which form an
open subspace in $\Bbb HP^n$, and on the space $\Cal N$ of isotropic
lines, which is a closed subspace of $\Bbb HP^n$. This identifies
$\Cal N$ with $G/P$, where $P\subset G$ is the stabilizer of an
isotropic line in $\Bbb H^{n+1}$.  

For a quaternionic line $\ell\subset\Bbb H^{n+1}$, we can realize
$T_\ell\Bbb HP^n$ as $\Bbb H^{n+1}/\ell$. Differentiating the defining
equation $\tau(v,v)=0$, we see that $\Cal N\subset\Bbb HP^n$ is a
smooth submanifold and that $T_\ell\Cal N=\{w\in
V:\im(\tau(v,w))=0\}/\ell$, so in particular this submanifold has real
codimension three.  Moreover, there is a natural subspace
$H_\ell:=\{w:\tau(v,w)=0\}/\ell\subset T_\ell\Cal N$, which has real
codimension one. Taking an explicit realization of $G$, it is easy to
verify that the codimension--one subbundle $H\subset T\Cal N$
constructed above defines a parabolic contact structure, consisting of
a quaternionic structure on $H$ such that the Levi--bracket is the
real part of a quaternionically skew Hermitian form. One further
verifies directly that the resulting subgroup $G_0$ is a two--fold
covering of the group $CSO^*(2n)$ generated by $SO^*(2n)$ and real
multiples of the identity.

On the other hand, consider the complex Grassmannian $Gr(2,\Bbb
C^{n+1})$ of two planes as a homogeneous space of $SU(n+1)$. Then it
is well known that this space admits an invariant complex structure as
well as an invariant quaternionic structure, and an invariant
Riemannian metric which is K\"ahler respectively quaternion--K\"ahler
with respect to these structures. In Corollary 3.8 of \cite{PCS1} we
have shown that the K\"ahler form of the K\"ahler metric together with
the quaternionic structure defines an $SU(n+1)$--invariant
PCS--structure of quaternionic type on $Gr(2,\Bbb C^{n+1})$.

\begin{thm}\label{thm2.6}
  Let $\Cal N\subset\Bbb HP^n$ be the space of quaternionic lines
  which are isotropic for a quaternionically skew Hermitian form. Then
  there is a projection $q:\Cal N\to Gr(2,\Bbb C^{n+1})$, which
  defines a global parabolic contactification of the
  $SU(n+1)$--invariant PCS--structure of type $(PSO^*(2n),P)$ on the
  Grassmannian with circles as fibers.
\end{thm}
\begin{proof}
  Let $\langle\ ,\ \rangle$ be the standard positive definite
  quaternionically Hermitian form on $V:=\Bbb H^{n+1}$. Fix a
  quaternionically linear map $\Cal J:V\to V$ such that $\Cal J\o\Cal
  J=-\id$ and such that $\langle \Cal J(v),\Cal J(w)\rangle=\langle
  v,w\rangle$. (For example, one can take multiplication by $i$ from
  the left.) Then it follows immediately that $\tau(v,w):=\langle
  v,\Cal J(w)\rangle$ is a non--degenerate quaternionically skew
  Hermitian form, so we can use this form to realize $SO^*(2n+2)$. A
  quaternionic line $\ell\subset V$ is isotropic for $\tau$ if and
  only if $\Cal J(\ell)\subset \ell^\perp$, the orthocomplement of $\ell$
  with respect to $\langle\ ,\ \rangle$. 

Scalar multiplication with respect to the complex structure on $V$
defined by $\Cal J$ defines an action of $U(1)$ on $V$ by
quaternionically linear maps. From the construction of $\tau$, one
readily verifies that this action is orthogonal for $\tau$, so we have
found a subgroup of $G$ isomorphic to $U(1)$. By definition,
$\tau(v,v)=0$ means that $\langle v,\Cal J(v)\rangle=0$, so $\Cal
J(v)$ does not lie in the quaternionic span of $v$. This shows that
only $\pm 1\in U(1)$ act as the identity on $G/P$. Moreover, for
$v\in\ell$ we have $\tau(v,\Cal J(v))=-\langle v,v\rangle\neq 0$, so
the infinitesimal automorphism generating this subgroup is transversal
on all of $G/P$.

Given an isotropic line $\ell\in\Cal N$, we conclude from $\ell\cap
\Cal J(\ell)=\{0\}$ that $\ell\oplus\Cal J(\ell)$ is quaternionic
subspace in $V$ of quaternionic dimension two, which in addition is
invariant under $\Cal J$. The $U(1)$--orbit of $\ell$ by construction
consists of null lines which are contained in this plane. From above,
we know that the space of null lines in a quaternionic plane is a
codimension three subspace in $\Bbb HP^1\cong S^4$, so it has
dimension one, and it is easily verified to be connected. This shows
that the $U(1)$--orbit of $\ell$ coincides with the space of those
isotropic lines which are contained in $\ell\oplus\Cal J(\ell)$.

Hence we see that mapping $\ell$ to $\ell\oplus\Cal J(\ell)$ defines a
(evidently smooth) map $q$ from $G/P$ to the space of $\Cal
J$--invariant quaternionic planes in $V$, whose fibers are exactly the
orbits of the $U(1)$--action constructed above. On the other hand,
since $\tau(v,\Cal J(v))\neq 0$ for any $v\in V$, it follows that the
restriction of $\tau$ to any $\Cal J$--invariant quaternionic plane
$W\subset V$ is non--degenerate. In particular, any such plane
contains a family of quaternionic null--lines isomorphic to $U(1)$, so
$q$ is surjective. 

So it remains to show that the space of $\Cal J$--invariant
quaternionic planes in $V$ is isomorphic to the complex
Grassmannian. Choose two anti--commuting imaginary unit quaternions,
for example $i$ and $j$. Then right multiplication by $i$ defines a
linear map $V\to V$ which squares to $-\id$ and commutes with $\Cal
J$. Hence we obtain a decomposition $V=V^{(1,0)}\oplus V^{(0,1)}$ into
two $\Cal J$--invariant summands characterized by $v\cdot i=\Cal J(v)$
for $v\in V^{(1,0)}$ and $v\cdot i=-\Cal J(v)$ for $v\in
V^{(0,1)}$. Moreover, one immediately verifies that multiplication by
$j$ from the right maps $V^{(1,0)}$ to $V^{(0,1)}$ and vice versa, so
both spaces must have complex dimension $n+1$. The projections onto
$V^{(1,0)}$ and $V^{(0,1)}$ can be constructed from $\Cal J$ and right
multiplication by $i$. A $\Cal J$--invariant quaternionic subspace
$W\subset V$ thus is invariant under the projections, and hence
$W=(W\cap V^{(1,0)})\oplus (W\cap V^{(1,0)})\cdot j$. Hence we see
that we can identify the space of $\Cal J$--invariant quaternionic
planes in $V$ with the space of complex planes in $V^{(1,0)}$ and
hence with the Grassmannian $Gr(2,\Bbb C^{n+1})$.

Hence we see that there is a natural projection $q:\Cal N\to Gr(2,\Bbb
C^{n+1})$, which is a quotient by a transversal infinitesimal
contactomorphism. The action of $G$ on $\Cal N$ does not descend to
the Grassmannian, but one can consider the stabilizer of $\Cal J$ in
$G$. It is easy to see that this stabilizer is isomorphic to $U(n+1)$ via
the action on $V^{(1,0)}$, and the action on $Gr(2,\Bbb C^{n+1})$
further factorizes to $SU(n+1)$. The conformally symplectic structure
on the Grassmannian is $SU(n+1)$ invariant and coincides with the one
used in Example 3.7 and Corollary 3.8 of \cite{PCS1}.

Now as in the first example, the projection $G\to G/P$ induces a
covering on each $U(1)$--orbit but this time, the covering is just
$2$--fold. Hence $q$ can be made into a PCS--quotient by passing to
the quotient of $G$ by $\Bbb Z_2=\{\pm\Bbb I\}$, and the resulting
PCS--structure of quaternionic type on the Grassmannian is exactly the
one from Corollary 3.8 of \cite{PCS1}.
\end{proof}

\subsection{Contactification of PCS--structures}\label{2.7}
For the applications to BGG sequences, we will mainly need the
following counterpart to Theorem \ref{thm2.5}. We show that any local
contactification of the conformally symplectic structure underlying a
PCS--structure can be canonically made into a PCS--quotient. Together
with the result in Lemma 3.1 of \cite{Cap-Salac}, this shows that any
PCS--structure can be locally realized as a PCS--quotient.

\begin{thm}\label{thm2.7}
  Suppose that $(M^\#,H)$ is a contact manifold, $(M,\ell)$ is a
  conformally symplectic manifold and $q:M^\#\to M$ is a reduction by
  a transverse infinitesimal contactomorphism $\xi\in\frak X(M^\#)$.

  Then any PCS--structure on $M$, which has $\ell$ as its underlying
  conformally symplectic structure canonically lifts to a parabolic
  contact structure on $M^\#$ for which $\xi$ is an infinitesimal
  automorphism.
\end{thm}
\begin{proof}
  As discussed in Section \ref{2.2}, a PCS--structure on $M$
  corresponding to the group $G_0\subset\ CSp(\frak g_{-1})$ is given
  by a principal $G_0$--bundle $p:\Cal G_0\to M$ together with a
  one--form $\th\in\Om^1(\Cal G_0,\frak g_{-1})$ which is
  $G_0$--equivariant and strictly horizontal. Now we can simply form
  the pullback $\Cal G_0^\#:=q^*\Cal G_0\to M^\#$, which is a
  principal $G_0$--bundle over $M^\#$. Explicitly, $\Cal
  G_0^\#=\{(u,x^\#)\in\Cal G_0\x M^\#:p(u)=q(x^\#)\}$, and the
  principal right action on this bundle is given by the principal
  right action of $\Cal G_0$ acting on the first factor. Pulling back
  $\th$ along the projection onto the first factor, we obtain a smooth
  $\frak g_{-1}$--valued one--form $\th^\#_{-1}$ on $\Cal G_0^\#$,
  which evidently is $G_0$--equivariant. From the construction it is
  clear that, for each $x^\#\in M^\#$, $\th_{-1}^\#(u,x^\#)$ descends
  to a linear map $T_{x^\#}M^\#\to\frak g_{-1}$, which restricts to a
  linear isomorphism on the contact subspace $H_{x^\#}$.

  Next, there is a unique contact form $\al$ on $M^\#$ such that
  $\al(\xi)\equiv 1$ and by Proposition 2.2 of \cite{Cap-Salac} there
  is a unique symplectic form $\om$ on $M$ which is a section of
  $\ell$ and satisfies $q^*\om=d\al$. By assumption, the conformally
  symplectic structure induced by the reduction $p:\Cal G_0\to M$ of
  structure group is $\ell\subset\La^2T^*M$. Given a point $u\in\Cal
  G_0$ with $p(u)=x$, we have $\om(x)\in\ell_x$ so there is a linear
  isomorphism $\ps(u):\frak g_{-2}\to\Bbb R$ such that, viewing
  $\th(u)$ as a map $T_xM\to\frak g_{-1}$, we get
$$
\om(x)(\eta_1,\eta_2)=\ps(u)([\th(u)(\eta_1),\th(u)(\eta_2)])
\quad\forall \eta_1,\eta_2\in T_xM,
$$ where the Lie bracket is in $\frak g_{-}$. This defines a smooth
map $\ps:\Cal G_0\to L(\frak g_{-2},\Bbb R)$ and equivariancy of $\th$
implies that $\th(u\cdot g)(\eta)=\Ad(g^{-1})(\th(u)(\eta))$, and
hence $\ps(u\cdot g)=\ps(u)\o\Ad(g)$.

Using this, we now define $\th_{-2}^\#\in\Om^1(\Cal G_0^\#,\frak
g_{-2})$ by $\th_{-2}^\#(u,x^\#)=\ps(u)^{-1}\o (p^\#_0)^*\al$, where
$p^\#_0:\Cal G_0^\#\to M^\#$ is the canonical projection. This evidently
is a smooth one--form, and since $(p^\#_0)^*\al$ is invariant under the
principal right action, we see that $\th^\#_{-2}$ is
$G_0$--equivariant. Next, the kernel of $\th^\#_{-2}$ in a point
coincides with the kernel of $(p^\#_0)^*\al$ in that point and thus with
the pre--image $(Tp^\#_0)^{-1}(H)$. Finally, if $\tilde\eta_1$ and
$\tilde\eta_2$ are two sections of this subbundle, then we compute
\begin{align*}
d\th_{-2}^\#(u,x^\#)(\tilde\eta_1,\tilde\eta_2)&=
-\th_{-2}^\#(u,x^\#)([\tilde\eta_1,\tilde\eta_2])=
-\ps(u)^{-1}((p^\#_0)^*\al(u,x^\#)([\tilde\eta_1,\tilde\eta_2])\\
&=\ps(u)^{-1}(d(p^\#_0)^*\al(u,x^\#)(\tilde\eta_1,\tilde\eta_2)). 
\end{align*}
Now since $d(p^\#_0)^*\al=(p^\#_0)^*d\al=(p^\#_0)^*q^*\om$, we get, putting
$x=q(x^\#)$ and denoting by $\eta_i\in T_xM$ the value of
$\tilde\eta_i(u,x^\#)$ under the natural
projection for $i=1,2$: 
\begin{align*}
  d\th_{-2}^\#(u,x^\#)(\tilde\eta_1,\tilde\eta_2)&=
  \ps(u)^{-1}(\om(x)(\eta_1,\eta_2))=
  [\th(u)(\eta_1),\th(u)(\eta_2)]\\
  &=[\th^\#_{-1}(u,x^\#)(\tilde\eta_1),\th^\#_{-1}(u,x^\#)(\tilde\eta_2)],
\end{align*}
where the Lie bracket is in $\frak g_-$. Together with the above
properties, this shows that $\th^\#:=(\th^\#_{-2},\th^\#_{-1})$ makes
$p^\#_0:\Cal G_0^\#\to M^\#$ into a regular infinitesimal flag structure
on $M^\#$, see Section \ref{2.3}. Since $\frak g$ is not of type
$C_n$, this is equivalent to a parabolic contact structure of type
$(G,P)$.

So it remains to show that $\xi$ is an infinitesimal automorphism of
this parabolic contact structure or equivalently of the infinitesimal
flag structure. Observe first that $T_{(u,x^\#)}\Cal
G_0^\#=\{(\eta_1,\eta_2)\in T_u\Cal G_0\x
T_{x^\#}M^\#:T_up\cdot\eta_1=T_{x^\#}q\cdot\eta_2\}$, and hence
$\tilde\xi:=(0,\xi)$ is a well--defined vector field on $\Cal
G_0^\#$. The flow of this vector field is clearly given by
$\Fl^{\tilde\xi}_t(u,x^\#)=(u,\Fl^\xi_t(x^\#))$ and it suffices to
prove that this flow preserves $\th^\#$.

Take a point $(u,x^\#)\in\Cal G_0^\#$ and a tangent vector
$(\eta_1,\eta_2)\in T_{(u,x^\#)}\Cal G_0^\#$. Then
$$
T_{(u,x^\#)}\Fl^{\tilde\xi}_t\cdot
(\eta_1,\eta_2)=(\eta_1,T_{x^\#}\Fl^\xi_t\cdot\eta_2),
$$
wherever the right hand side is defined. But by definition 
$$
\th^\#_{-1}(u,x^\#)(\eta_1,\eta_2)=\th(u)(\eta_1), 
$$ 
which immediately implies
$(\Fl^{\tilde\xi}_t)^*\th^\#_{-1}=\th^\#_{-1}$, whenever the flow is
defined. Again by definition, we get
$$
\th_{-2}^\#(u,x^\#)(\eta_1,\eta_2)=\ps(u)^{-1}(\al(x^\#)(\eta_2)),
$$
and hence $(\Fl^{\tilde\xi}_t)^*\th^\#_{-2}=\th^\#_{-2}$ follows
immediately from $(\Fl^{\xi}_t)^*\al=\al$. 
\end{proof}

\subsection{Local uniqueness of PCS--contactifications}\label{2.8} 
To complete the picture, we have to prove that different realizations
of PCS--structures as PCS--quotients as constructed in Theorem
\ref{thm2.7} are locally compatible. The corresponding result for
conformally symplectic structures has been proved in Proposition 3.1
of \cite{Cap-Salac}, so we only have to prove compatibility with the
additional structures.

\begin{thm}\label{thm2.8} 
  Suppose that $M$ and $\tilde M$ carry a PCS--structure of some fixed
  type with underlying conformally symplectic structures $\ell$ and
  $\tilde\ell$. Suppose further that $(M^\#,H)$ and $(\tilde
  M^\#,\tilde H)$ are contact manifolds and that $q:M^\#\to M$ and
  $\tilde q:\tilde M^\#\to\tilde M$ are reductions with respect to
  transverse symmetries $\xi\in\frak X(M^\#)$ and $\tilde\xi\in\frak
  X(\tilde M^\#)$. Suppose finally that $\ph:M\to\tilde M$ is a
  PCS--diffeomorphism and that $\ph^\#:M^\#\to\tilde M^\#$ is a
  contactomorphism which lifts $\ph$. 

Then $(\ph^\#)^*\tilde\xi=\la\xi$ for a nowhere--vanishing, locally
constant function $\la$ and $\ph^\#$ is an automorphism of the lifted
parabolic contact structures from Theorem \ref{thm2.7}.
\end{thm}
\begin{proof}
Since $\ph^\#$ is a contactomorphism, $(\ph^\#)^*\tilde\xi\in\frak
X(M^\#)$ is a transverse infinitesimal automorphism of the contact
structure $H\subset TM^\#$. On the other hand, since $\ph^\#$ lifts $\ph$ it
has to map fibers of $q$ to fibers of $\tilde q$, which shows that
$(\ph^\#)^*\tilde\xi=\la\xi$ for some smooth function $\la$ on
$M^\#$. But then for $\eta\in\frak X(M^\#)$ we get $[\la
  \xi,\eta]=\la[\xi,\eta]-(\eta\cdot\la)\xi$, and if $\eta\in\Ga(H)$,
then by assumption both the left hand side and the first summand in
the right hand side are sections of $H$. Since $\xi$ is transverse,
this implies that $\eta\cdot \la=0$ for any $\eta\in\Ga(H)$ and since
$H$ is bracket--generating this shows that $\la$ is locally constant. 

Now let $(p:\Cal G_0\to M,\th)$ and $(\tilde p:\tcg_0\to\tilde
M,\tilde\th)$ be the bundles defining the PCS--structures. Then by
assumption there is an isomorphism $\Ph:\Cal G_0\to\tcg_0$ of
principal $G_0$--bundles such that $\tilde p\o\Ph=\ph\o p$ and such
that $\Ph^*\tilde\th=\th$. Forming the pullback bundles as in the
proof of Theorem \ref{thm2.7}, we see that
$$
(u,x^\#)\mapsto (\Ph(u),\ph^\#(x^\#))
$$ 
defines an isomorphism $(\Ph,\ph^\#):\Cal G_0^\#\to\tcg_0^\#$ lifting
$\ph^\#$, so it remains to prove that this is compatible with the
frame forms on the two bundles constructed in the proof of Theorem
\ref{thm2.7}. Inserting the definitions, one immediately verifies that
$$ 
(\Ph,\ph^\#)^*\tilde\th_{-1}^\#(u,x^\#)(\eta_1,\eta_2)=
\Ph^*\tilde\th(u)(\eta_1)=\th(u)(\eta_1),
$$ 
and thus $(\Ph,\ph^\#)^*\tilde\th_{-1}^\#=\th_{-1}^\#$. On the other
hand, since $(\ph^\#)^*\tilde\xi=\la\xi$, we see that the contact
forms corresponding to $\xi$ and $\tilde\xi$ are related by
$(\ph^\#)^*\tilde\al=\frac1{\la}\al$. Since $\frac1\la$ is locally
constant, the corresponding symplectic forms are related by
$\ph^*\tilde\om=\frac1{\la}\om$. Viewing the values of $\th$ and
$\tilde\th$ as linear isomorphisms on tangent spaces of $M$ and
$\tilde M$, the fact that $\Ph^*\tilde\th=\th$ reads as
$$
\tilde\th(\Ph(u))(T_{p(u)}\ph\cdot\eta)=\th(u)(\eta)
$$ 
for all $u\in\Cal G_0$ and $\eta\in T_{p(u)}M$. Together with
$\ph^*\tilde\om=\frac1{\la}\om$, this implies that the isomorphisms
$\ps$ and $\tilde\ps$ as constructed in the proof of Theorem \ref{thm2.7}
satisfy $\tilde \ps(\Ph(u))=\tfrac1\la\ps(u)$. Using this, we compute 
\begin{align*}
  (\Ph,\ph^\#)^*\tilde\th_{-2}^\#&(u,x^\#)(\eta_1,\eta_2)=
  \tilde\th_{-2}^\#(\Ph(u),\ph^\#(x^\#))(T_u\Ph\cdot\eta_1,
  T_{x^\#}\ph^\#\cdot\eta_2)\\ &=\tilde\ps(\Ph(u))^{-1}((\ph^\#)^*\tilde\al(x^\#)(\eta_2))=\la\ps(u)^{-1}(\tfrac1\la\al(x^\#)(\eta_2)),
\end{align*}
so $(\Ph,\ph^\#)^*\tilde\th_{-2}^\#=\th_{-2}^\#$.
\end{proof}

\subsection{Distinguished connections and
  contactifications}\label{2.9}  
As shown in \cite{PCS1}, any PACS--structure on $M$ determines a
canonical compatible linear connection on $TM$. In the case of
PCS--structures, these are closely related to special symplectic
connections as studied in \cite{Cahen-Schwachhoefer}, see Section 4.7
of \cite{PCS1} for a discussion of the relation. We now conclude the
discussion of parabolic contactifications by relating the
distinguished connection associated to a PCS--quotient of a parabolic
contact structure to distinguished connections ``upstairs''. On the
one hand, this provides alternative proofs for and extensions of some
results from \cite{Cahen-Schwachhoefer}, on the other hand, the result
will be useful for the discussion of conformally Fedosov structures
in Section \ref{3} below as well as for the third part of this series.

A family of distinguished connections associated to a parabolic
contact structure is obtained via so--called Weyl structures as
discussed in Chapter 5 of \cite{book}. Consider the Cartan geometry
$(\Cal G^\#\to M^\#,\om)$ of type $(G,P)$ determined by a parabolic
contact structure and the underlying bundle $\Cal G^\#_0\to M$. A
\textit{Weyl structure} is defined as a $G_0$--equivariant section
$\si$ of the natural projection $\Cal G^\#\to\Cal G_0^\#$. Taking the
component $\om_0\in\Om^1(\Cal G^\#,\frak g_0)$ of the Cartan
connection in $\frak g_0\subset\frak g$, the pull back $\si^*\om_0$ is
a principal connection on $\Cal G^\#_0$ called the \textit{Weyl
  connection} determined by $\si$.

A principal connection on $\Cal G_0^\#$ can be equivalently viewed as
a linear connection $\nabla^\#$ on the vector bundle $H\to M^\#$,
which is compatible with the reduction of structure group defined by
$\Cal G_0^\#$. In particular, $\nabla^\#$ is a \textit{contact
  connection}, in the sense that the induced connection on $\La^2H$
preserves the subbundle $\La^2_0H$ formed by the kernel of the
Levi--bracket.

Further, the pull back of the components in $\frak g_{-2}\oplus\frak
g_{-1}$ of the Cartan connection $\om$ along $\si$ can be interpreted
as defining an isomorphism between $TM^\#$ and its associated graded
vector bundle $H\oplus Q$, where $Q:=TM^\#/H$. In particular, this
defines a projection $\pi$ from $TM^\#$ onto the subbundle $H\subset
TM^\#$. Using this projection, we can define the \textit{contact
  torsion} of $\nabla^\#$ via
$$
T(\eta_1,\eta_2)=\nabla^\#_{\eta_1}\eta_2-\nabla^\#_{\eta_2}\eta_1-
\pi([\eta_1,\eta_2])
$$
for sections $\eta_1,\eta_2\in\Ga(H)$. This expression is immediately
seen to be bilinear over smooth functions and thus defines a bundle
map $T:\La^2H\to H$.

We have to recall a few further facts on Weyl structures from Section
5.2.11 of \cite{book}. The line bundle $Q$ can be used as a so--called
bundle of scales, so passing to the induced linear connection on $Q$
defines a bijection between Weyl structures and linear connections on
$Q$. Now consider a reduction $q:M^\#\to M$ by a transverse symmetry
$\xi\in\frak X(M^\#)$. Then we can project $\xi$ to a section of $Q$,
which by construction is nowhere vanishing. Of course, there is a
unique linear connection on $Q$, for which this section is parallel,
and this in turn uniquely determines a Weyl structure $\si$. This will
be called the \textit{Weyl structure associated to $\xi$}.

\begin{thm}\label{thm2.9}
  Let $q_0:\Cal G_0^\#\to\Cal G_0$ be a PCS--quotient with respect to
  a transverse infinitesimal automorphism $\xi_0\in\frak X(\Cal
  G_0^\#)$ with base map $q:M^\#\to M$.

  Then there is a unique principal connection form $\ga^\#$ on
  $p^\#_0:\Cal G_0^\#\to M^\#$ for which $\xi_0\in\frak X(\Cal
  G_0^\#)$ is horizontal and which coincides with the principal
  connection form determined by the Weyl structure associated to $\xi$
  on $T^{-1}\Cal G_0^\#\subset T\Cal G_0^\#$. The form $\ga^\#$
  descends to a principal connection form on $\Cal G_0\to M$, which
  corresponds to the canonical connection associated to the PCS
  structure from Theorem \ref{thm2.2}.
\end{thm}
\begin{proof}
  Let $\si:\Cal G^\#_0\to\Cal G^\#$ be the Weyl structure determined
  by $\xi$ and let $\tilde\ga^\#=\si^*\om_0$ be the corresponding Weyl
  connection. Since $\xi_0$ is $G_0$--invariant, the function
  $\tilde\ga^\#(\xi_0):\Cal G_0^\#\to\frak g_0$ is
  $G_0$--equivariant. Next, let $\al\in\Om^1(M^\#)$ be the unique
  contact form with $\al(\xi)=1$ and put
  $\ga^\#:=\tilde\ga^\#-\tilde\ga^\#(\xi_0)(p_0^\#)^*\al$. This is
  evidently $G_0$--equivariant, and since $(p_0^\#)^*\al$ vanishes
  on $T^{-1}\Cal G_0^\#$, it is a principal connection form which
  coincides with $\tilde\ga^\#$ on $T^{-1}\Cal G_0^\#$. Since $\xi_0$
  lifts $\xi$, we get $\ga^\#(\xi_0)=0$, so we have verified the
  required properties. Of course, these pin down $\ga^\#$ uniquely.

  By definition, the fibers of $q_0:\Cal G_0^\#\to\Cal G_0$ are
  connected and can locally be represented as flow lines of
  $\xi_0$. Hence by Corollary 2.3 in \cite{BCG3} we can prove that
  $\ga^\#$ descends to $\Cal G_0$ by showing that it is invariant
  under the flow of $\xi_0$. As discussed in Section \ref{2.3}, the
  infinitesimal automorphism $\xi_0$ uniquely lifts to a
  $P$--invariant vector field $\tilde\xi\in\frak X(\Cal G)$ whose
  local flows preserve the Cartan connection $\om$. Supposing that
  both flows are defined, the map
  $\Fl^{\tilde\xi}_{-t}\o\tilde\si\o\Fl^{\xi_0}_t$ is again a Weyl
  structure. Since $(\Fl^{\tilde\xi}_t)^*\om=\om$, we conclude that
  $(\Fl^{\xi_0}_t)^*\tilde\ga^\#$ is the Weyl connection associated to
  this pulled back Weyl structure. Thus we can prove invariance of
  $\tilde\ga^\#$ under the flow of $\xi_0$ by showing that
  $\Fl^{\tilde\xi}_{-t}\o\tilde\si\o\Fl^{\xi_0}_t=\si$ for small $t$.

  Denoting by $\tau\in\Ga(Q)$ the section obtained by projecting
  $\xi$, it is evident that $(\Fl^\xi_t)^*\tau$ is parallel for the
  Weyl connection corresponding to
  $\Fl^{\tilde\xi}_{-t}\o\tilde\si\o\Fl^{\xi_0}_t$, and we know from
  above that this characterizes the pulled back Weyl structure. But
  since $\tau$ is a projection of $\xi$ and the Lie derivative of
  $\xi$ along $\xi$ equals $[\xi,\xi]=0$, we conclude
  $(\Fl^\xi_t)^*\tau=\tau$, so
  $\Fl^{\tilde\xi}_{-t}\o\tilde\si\o\Fl^{\xi_0}_t=\si$.

  Knowing that $\tilde\ga^\#$ is invariant under the flow of $\xi_0$,
  it follows readily that the same is true for the function
  $\tilde\ga^\#(\xi_0)$. Also, the pullback of $\al$ clearly is
  invariant under the flow of $\xi_0$, since $\al$ is invariant under
  the flow of $\xi$. Thus we conclude that $\ga^\#$ is invariant under
  the flow of $\xi_0$ and hence descends to a $\frak g_0$--valued
  one--form, which then clearly is a principal connection form on
  $\Cal G_0\to M$.

  To complete the proof it thus suffices to prove that the
  corresponding linear connection on $TM$ has algebraically harmonic
  torsion. Now the torsion of a linear connection on $TM$ can be
  computed from the corresponding principal connection by evaluating
  the exterior derivative of the soldering form on horizontal lifts of
  the vector fields. Now the soldering form on $\Cal G_0$ pulls back
  to the form $\th_{-1}^\#$ on $\Cal G_0^\#$. This can be extended to
  a one--form on $\Cal G_0^\#$ by requiring that it vanishes on
  $\xi_0$ (which corresponds to the isomorphism between $TM^\#$ and
  $H\oplus Q$ defined by the Weyl structure determined by
  $\xi$). Evaluating the exterior derivative of this one--form on the
  horizontal lifts sections of $H\to M^\#$, one obtains the function
  corresponding to the contact torsion of the corresponding Weyl
  connection. This contact torsion is well known to have values in
  $\ker(\square)$, compare with Theorem 5.2.11 of \cite{book}. 

  Now consider the definition of the contact torsion $T(\eta_1,\eta_2)$
  in case that $\eta_1$ and $\eta_2$ are local lifts of vector fields
  on $M$. Then $[\eta_1,\eta_2]$ is a local lift of their Lie bracket
  and this differs from $\pi([\eta_1,\eta_2])$ by a multiple of $\xi$.
  Thus $\pi([\eta_1,\eta_2])$ is the section of $H$ lifting the Lie
  bracket of the downstairs vector fields. This immediately implies
  that this contact torsion coincides with the torsion of the linear
  connection associated to $\ga$, which completes the proof.
\end{proof}

\subsection{Contactifications and special symplectic
  connections}\label{2.10} 
Using Theorem \ref{thm2.9}, we can now complete the discussion of the
relation between PCS--structures and special symplectic connections in
the sense of \cite{Cahen-Schwachhoefer}. Consider a contact grading of
a simple Lie algebra $\frak g$, which is not of type $C_n$, and let
$\frak g_0\subset\mathfrak{csp}(\frak g_{-1})$ be the corresponding
reductive subalgebra. As discussed in Section 4.7 of \cite{PCS1}, the
intersection $\frak g^0_0$ of $\frak g_0$ with $\frak{sp}(\frak
g_{-1})$ is called a special symplectic subalgebra and there is a
corresponding concept of special symplectic connections. In particular,
that class contains all linear connections with exceptional holonomy
that preserve a symplectic form. 

In Theorem 4.7 of \cite{PCS1}, we have shown that any special
symplectic connection of type $\frak g$ is the distinguished
connection of some torsion--free PCS--structure of type $(G,P)$ for a
group $G$ with Lie algebra $\frak g$. Now we can characterize the
special symplectic connections among these distinguished connections
via local flatness of contactifications.

\begin{thm}\label{thm2.10}
  Consider a PCS--structure of type $(G,P)$ on a smooth manifold $M$
  and let $\nabla$ be the distinguished connection on $TM$ determined
  by this structure. Then $\nabla$ is a special symplectic connection
  if and only if any local parabolic contactification of $M$ is
  locally flat (as a parabolic contact structure).
\end{thm}
\begin{proof}
  For a local parabolic contactification $M^\#$ of $M$ with
  corresponding infinitesimal automorphism $\xi$, we have seen in the
  last part of the proof of Theorem \ref{thm2.9} that the contact torsion
  of the Weyl connection $\nabla^\#$ determined by $\xi$ coincides
  with the torsion of $\nabla$. Torsion--freeness of $\nabla$ thus
  implies vanishing of the contact torsion of $\nabla^\#$. But by
  Lemma 4.2.2 of \cite{book}, all the harmonic torsion of a parabolic
  contact structure is contained in the contact torsion, so
  torsion--freeness of $\nabla$ is equivalent to vanishing of the
  harmonic torsion of the parabolic geometry on $M^\#$.

  From the discussion of parabolic contact structures in Section 4.2
  of \cite{book} it follows that for $\frak g$ not of type $A_n$,
  vanishing of the harmonic torsion is equivalent to local flatness of
  the parabolic geometry. On the other hand, by Theorem 4.7 of
  \cite{PCS1}, for $\frak g$ not of type $A_n$, torsion--freeness of
  $\nabla$ is equivalent to $\nabla$ being special symplectic, so the
  proof is complete in this case.

  For $\frak g$ of type $A_n$, Proposition 4.2.3 and the discussion in
  Section 4.2.4 in \cite{book} show that vanishing of the harmonic
  torsion of the parabolic geometry on $M^\#$ is equivalent to
  torsion--freeness of this geometry. Moreover, we can apply part (3)
  of Theorem 4.2.2 of \cite{book} to $\nabla^\#$ and conclude that the
  harmonic curvature component in homogeneity two of this geometry
  (which is the only remaining component in dimensions $\geq 5$) can
  be recovered as a component of the curvature of $\nabla^\#$. This
  curvature is a two--form with values in $\text{End}(H)$ and one has
  to restrict this two--form to $\La^2_0H^*\subset\La^2 T^*M$, the
  kernel of the Levi--bracket within $\La^2H$. This defines a section
  of the bundle associated to $\La^2_0\frak g_1\otimes\frak g_0$, and
  one has to further restrict to the kernel of the Kostant--Laplacian
  $\square$ on this representation. But on $\La^2_0H^*$, one can also
  compute the curvature via
$$
R(\eta_1,\eta_2)(\eta_3)=\nabla^\#_{\eta_1}\nabla^\#_{\eta_2}\eta_3-
\nabla^\#_{\eta_2}\nabla^\#_{\eta_1}\eta_3-
\nabla^\#_{\pi([\eta_1,\eta_2])}\eta_3.
$$
This shows that the relevant component of the curvature of $\nabla^\#$
is just the pullback of the curvature of $\nabla$. Hence we see that
the parabolic contact structure on $M^\#$ is locally flat if and only
if the component of its curvature in the subbundle corresponding to
$\ker(\square)\subset \La^2\frak g_1\otimes\frak g_0$ vanishes
identically. This subspace is well known to be the $\frak
g_0$--irreducible component of maximal highest weight, which
immediately implies that this condition is equivalent to vanishing of
the Bochner curvature of $\nabla$ and thus the result.
\end{proof}

\section{Conformally Fedosov structures and their
  contactifications}\label{3} 

In the previous discussions, we have excluded type $C_n$, since in
this case $\frak g_0=\frak{csp}(\frak g_{-1})$. This corresponds to
the fact that the parabolic contact structure determined by this
grading is not defined by a regular infinitesimal flag structure in
the sense of Section \ref{2.3}. The contact grading of
$\frak{sp}(2n,\Bbb R)$ gives rise to so--called contact projective
structures, an analog of classical projective structures. In this
section, we discuss a geometric structure refining a conformally
symplectic structure, which can be thought of as the PCS--structure
associated to the contact grading of $\frak{sp}(2n,\Bbb R)$. (It seems
to be possible to define a similar concept related to almost
conformally symplectic structures, but since contactifications are
only available in the conformally symplectic case, we do not pursue
this direction here.)

Initially, this structure can be described as a projective structure
compatible with a conformally symplectic structure. However, there is
a counterpart to the vanishing of the first prolongation for
PACS--structures, which implies that such a geometry is given by a
conformally symplectic structure and a torsion--free affine connection
preserving this structure. Locally, such a structure is
even equivalent to a symplectic structure together with a
torsion--free symplectic connection. These geometries were originally
introduced by M.G.~Eastwood and J.~Slov\'ak (in a slightly different
presentation) in the first version of \cite{E-S} under the name
``conformally Fedosov structures'', and we will keep this name here.

\subsection{Projective structures and conformally Fedosov
  structures}\label{3.1}
Recall that two torsion--free linear connections $\nabla$ and
$\hat\nabla$ on the tangent bundle of a smooth manifold $M$ are called
\textit{projectively equivalent} if and only if they have the same
geodesics up to reparametrization. This can be equivalently
characterized as the existence of a one--form $\Ups$ on $M$ such that
$$
\hat\nabla_\xi\eta=\nabla_\xi\eta+\Ups(\xi)\eta+\Ups(\eta)\xi,
$$ 
for all $\xi,\eta\in\frak X(M)$.  A projective structure on a smooth
manifold $M$ is then defined as a projective equivalence class of
torsion--free connections on $M$. Equivalently, it can be viewed as
being given by the family of one--dimensional submanifolds formed by
the geodesic paths of the connections in the class.

\begin{definition}\label{def3.1}
  A \textit{conformally Fedosov structure} on a conformally symplectic
  manifold $(M,\ell\subset\La^2T^*M)$ is a projective structure
  $[\nabla]$ on $M$, which contains a connection that preserves $\ell$.
\end{definition}

The following characterization is essentially contained in \cite{E-S},
we reproduce it for convenience of the reader. We use abstract index
notation as in \cite{E-S}.

\begin{prop}\label{prop3.1}
  For a projective equivalence class $\Ph$ of torsion--free affine
  connections on a conformally symplectic manifold $(M,\ell)$ of
  dimension $\geq 4$, the following conditions are equivalent:

  (i) $\Ph$ defines a conformally Fedosov structure on $(M,\ell)$.

  (ii) For one (or equivalently any) connection $\nabla=\nabla_a$ in
  the class $\Ph$ and one non--vanishing (or equivalently any) local
  section $\om=\om_{bc}$ of $\ell$, we have
  $\nabla_{(a}\om_{b)c}=\ph_{(a}\om_{b)c}$ for some one--form
  $\ph=\ph_a$.

  (iii) There is a unique connection $\nabla_a$ in $\Ph$ which
  preserves $\ell$.

\smallskip

\noindent
Moreover, the connection on $\ell$ induced by the connection $\nabla$
from (iii) is flat with any local closed section being parallel. In
particular, locally $\nabla$ is a torsion--free symplectic connection
with respect to each of these closed sections.
\end{prop}
\begin{proof}
  The formula for the change of the induced connection on two--forms
  from \cite{E-S} says that changing from $\nabla_a$ to $\hat\nabla_a$
  corresponding to a one--form $\Ups_a$ one gets
$$
\hat\nabla_a\om_{bc}=\nabla_a\om_{bc}-2\Ups_a\om_{bc}-\Ups_b\om_{ac}-\Ups_c\om_{ba}.
$$ 
If we symmetrize over $(a,b)$,then the last term does not contribute
and we get
$$
\hat\nabla_{(a}\om_{b)c}=\nabla_{(a}\om_{b)c}-3\Ups_{(a}\om_{b)c}. 
$$
This shows that if the condition in (ii) is satisfied for one
connection in the projective class $\Ph$, then it is satisfied for all
such connections. Likewise, one easily verifies that if the condition
in (ii) is satisfied for one non--vanishing section, then it is
satisfied for any section of $\ell$.

Condition (i) by definition means that there is a connection
$\nabla_a$ in $\Ph$ such that for any local section $\om_{bc}$ of
$\ell$, we have $\nabla_a\om_{bc}=\ph_a\om_{bc}$ for some one--form
$\ph=\ph_a$. Thus (i) implies (ii) and evidently (iii) implies (i), so
we can complete the proof by showing that (ii) implies (iii). 

To do this, take a local non--vanishing section $\om_{bc}$ of $\ell$,
which is closed as a two--form on $M$, and any connection $\nabla_a$
in the class $\Ph$. Then by (ii) there is a one--form $\ph_a$ such
that $\nabla_{(a}\om_{b)c}=\ph_{(a}\om_{b)c}$. Passing to the
projectively equivalent connection $\hat\nabla$ corresponding to
$\Ups_a=\tfrac13\ph_a$ we see from above that
$\hat\nabla_{(a}\om_{b)c}=0$. This means that apart from the apparent
skew symmetry in $b$ and $c$, $\hat\nabla_a\om_{bc}$ is also skew
symmetric in $a$ and $b$. But this easily implies that the complete
alternation of $\hat\nabla_a\om_{bc}$ is a multiple of
$\hat\nabla_a\om_{bc}$, and since this alternation coincides with
$d\om=0$, we see that $\om$ is parallel for $\hat\nabla$. This implies
the existence part of (iii) as well as the last claim.

To prove the uniqueness part of (iii), observe that by part (4) of
Proposition 2.3 of \cite{PCS1} every torsion--free connection
compatible with $\ell$ has the property that any local closed section
of $\ell$ is parallel for the induced connection. But if $\om$ is any
such section, then by non--degeneracy $\om^n=\om\wedge\dots\wedge\om$
is a volume form on $M$, which clearly is parallel for the induced
connection. But it is a well known fact in the theory of projective
structures that a projective class can contain at most one connection
which leaves some fixed volume form parallel.
\end{proof}

\subsection{Contact projective structures}\label{3.2}
A parabolic contact structure corresponding to the contact grading of
a simple Lie algebra of type $C_n$ is a so--called \textit{contact
  projective structure}. As discussed above, this geometry is
exceptional among parabolic contact structures, since it is not
equivalent to an underlying regular infinitesimal flag
structure. Rather than that, it can be described in terms of an
equivalence class of compatible (partial) connections.

To recall the necessary notions, let $M^\#$ be a smooth manifold
endowed with a contact structure $H\subset TM^\#$. A \textit{partial
  connection} on $H$ is a bilinear operator
$\nabla:\Ga(H)\x\Ga(H)\to\Ga(H)$, which is linear over smooth
functions in the first variable and satisfies the usual Leibniz rule
in the second variable. So this is like a covariant derivative, but it
is possible to differentiate only in directions lying in the contact
subbundle. It is no problem to form the induced partial connection on
$\La^2H$, and parallel to Section \ref{2.9} a \textit{partial contact
  connection} is one that preserves $\La^2_0H$. Since the
Levi--bracket identifies $\La^2H/\La^2_0H$ with $Q=TM^\#/H$, we get an
induced partial connection $\nabla^Q$ on the line bundle $Q$.

Using this, one can next associate a projection $\pi:TM^\#\to H$ to a
partial contact connection $\nabla$, see Section 4.2.2 and 4.2.6 of
\cite{book}. This is uniquely characterized by the fact that for
$\zeta\in\frak X(M^\#)$ and $\eta\in\Ga(H)$ we have
$$
\Cal L(\pi(\zeta),\eta)=\nabla^Q_\eta(\zeta+H)-([\eta,\zeta]+H).
$$ 
here $\zeta+H$ denotes the section of $Q$ determined by $\zeta$ and
likewise for the Lie bracket. Having this projection we define the
contact torsion $T:\La^2H\to H$ of a partial contact connection
$\nabla$ exactly as in Section \ref{2.9}.

Finally, there is a notion of contact projective equivalence for
partial contact connections. Consider a smooth section $\Ups$ of the
bundle $H^*$ of linear functionals on the contact subbundle. Then
using the Levi--bracket $\Cal L$, we define $\tilde\Ups:Q\to H$ by
$\Cal L(\tilde\Ups(\tau),\eta)=\Ups(\eta)\tau$. Given a partial
contact connection $\nabla$ on $H$, one defines a partial connection
$\hat\nabla$ on $H$ by
$$
\hat\nabla_{\eta_1}\eta_2=\nabla_{\eta_1}\eta_2+\Ups(\eta_1)\eta_2+
\Ups(\eta_2)\eta_1+\tilde\Ups(\Cal L(\eta_1,\eta_2)). 
$$
One easily verifies directly that this is again a partial contact
connection on $H$. Further, one shows that the projection $\hat\pi$
associated to $\hat\nabla$ is given by
$\hat\pi(\zeta)=\pi(\zeta)+2\tilde\Ups(\zeta+H)$ and this easily
implies that $\hat\nabla$ and $\nabla$ have the same contact torsion.

One calls two partial contact connections \textit{contact projectively
  equivalent} if they are related in the above way. A (torsion--free)
contact projective structure on $M^\#$ is then given by a contact
structure together with a class of contact projectively equivalent
partial contact connections with vanishing contact torsion. It turns
out that such a structure can be equivalently described by a Cartan
geometry as discussed in Section \ref{2.3}, which corresponds to the
contact grading of the Lie algebra $\frak{sp}(2n+2,\Bbb R)$. As a
group, one can take $G=Sp(2n+2,\Bbb R)$ and $P\subset G$ the
stabilizer of a ray in the standard representation $\Bbb R^{2n+2}$ of
$G$. Then it turns out that $G_0$ is the conformally symplectic group
$CSp(2n,\Bbb R)$, see Section 4.2.6 of \cite{book}.

\begin{remark}\label{rem3.2}
  (1) Similarly to classical projective structures, there is also an
  interpretation of contact projective equivalence in terms of
  geodesics in directions tangent to the contact distribution, and a
  corresponding description of contact projective structures, see
  \cite{Fox} and Section 4.2.7 in \cite{book}.

  (2) In contrast to the case of usual affine connections, it is in
  general not possible to remove the contact torsion of a partial
  contact connection without changing its contact geodesics. Hence it
  is natural to extend the notion of a contact projective structure to
  the case of non--vanishing contact torsion, see \cite{Fox}. It turns
  out that also in this more general setting, one can associate a
  canonical Cartan geometry to such a structure. The resulting Cartan
  geometry does not fit into the general scheme of regular normal
  parabolic geometries, however.
\end{remark}

\subsection{Contactifications of conformally Fedosov
  structures}\label{3.3}
Our final task in this article is to establish analogs of Theorems
\ref{thm2.5}, \ref{thm2.7} and \ref{thm2.8} for contact projective
structures and conformally Fedosov structures. This turns out to be
rather easy since, as indicated in Section \ref{3.1}, the projective
freedom in a conformally Fedosov structure does not really show
up. Likewise, in the setting of contact projective structures, we will
always have a preferred infinitesimal automorphism around. As
discussed in Section \ref{2.9}, this determines a preferred Weyl
structure and thus a distinguished representative in the contact
projective class. Hence we can always related connections rather than
projective equivalence classes, and this works similarly as in the
proof of Theorem \ref{thm2.9}.

\begin{thm}\label{thm3.3}
  Let $q:M^\#\to M$ be a quotient of a contact manifold by a
  transverse infinitesimal contactomorphism $\xi\in\frak X(M^\#)$, and
  let $\ell\subset\La^2T^*M$ be the induced conformally symplectic
  structure.

  (1) A contact projective structure on $M^\#$, for which $\xi$ is an
  infinitesimal automorphism, induces a conformally Fedosov structure
  on $M$ with underlying conformally symplectic structure $\ell$.

  (2) Conversely, a conformally Fedosov structure on $M$ with
  underlying conformally symplectic structure $\ell$ canonically lifts
  to a contact projective structure for which $\xi$ is an
  infinitesimal automorphism.

  (3) Suppose that we have two conformally Fedosov structures realized
  as quotients of contact projective structures as in (2). Then any
  local lift of a diffeomorphism respecting the conformally Fedosov
  structures to a contactomorphism is automatically compatible with
  the infinitesimal automorphisms as in Theorem \ref{thm2.8} and a local
  isomorphism of contact projective structures.
\end{thm}
\begin{proof}
  Let $\Cal G_0^\#\to M^\#$ be the natural frame bundle for the
  contact subbundle $H$ with structure group $G_0=CSp(2n,\Bbb R)$ and
  let $\Cal G_0\to M$ be the natural frame bundle for $TM$ induced by
  $\ell$, which also has structure group $G_0$, see Section 2.3 of
  \cite{PCS1}. By construction, there is a tautological soldering form
  $\th\in\Om^1(\Cal G_0,\Bbb R^{2n})$. Denoting as before by
  $T^{-1}\Cal G_0^\#$ the preimage of $H\subset TM^\#$ in $T\Cal
  G_0^\#$, we likewise get a tautological soldering form
  $\th_{-1}^\#\in\Ga(L(T^{-1}\Cal G_0^\#,\Bbb R^{2n}))$.

  By definition, the local flows of $\xi$ are contactomorphisms and
  hence lift to local principal bundle automorphisms of $\Cal G_0^\#$.
  Differentiating these local flows, we obtain a $G_0$--invariant lift
  $\xi_0\in\frak X(\Cal G_0^\#)$ of $\xi$. As observed in the proof of
  Theorem \ref{thm2.5} the tangent maps of $q$ restrict to linear
  isomorphisms on the contact subspaces, which provides a
  $G_0$--equivariant lift $q_0:\Cal G_0^\#\to\Cal G_0$ of $q:M^\#\to
  M$. Moreover, by construction the restriction of $q_0$ to each fiber
  of $\Cal G_0^\#$ is injective. As in the proofs of Lemma \ref{lem2.4}
  and Theorem \ref{thm2.5}, this shows that $q_0$ is a surjective
  submersion and that the restriction of $p_0^\#:\Cal G_0^\#\to M^\#$
  to any fiber of $q_0$ is a diffeomorphism onto a fiber of $q$, so in
  particular $q_0$ has connected fibers. Finally, it is easy to see
  that the tangent spaces to these fibers are spanned by $\xi_0$. 

  (1) We can pass to the Cartan geometry determined by the given
  contact projective structure on $M^\#$, which has $\Cal G_0^\#$ as
  its underlying $G_0$--bundle. The general theory of Weyl structures
  applies to contact projective structures. Hence as in the proof of
  Theorem \ref{thm2.9}, we see that the Weyl structure determined by
  $\xi$ determines a principal connection form $\tilde\ga^\#$ on $\Cal
  G_0^\#$. Here the Weyl connection is an extension of a partial
  contact connection with vanishing contact torsion, which lies in the
  contact projective class. Continuing as in the proof of Theorem
  \ref{thm2.9}, we can change this to a principal connection form
  $\ga^\#$ for which $\xi_0$ is horizontal, and show that this
  descends to $\Cal G_0$. Still as in that proof, we can verify that
  the resulting connection on $M$ is torsion--free and by construction
  it preserves $\ell$. So we can use the projective class of this
  connection to define a conformally Fedosov structure on $M$.

  (2) Starting with a conformally Fedosov structure on $M$, we know from
  Proposition \ref{prop3.1} that there is a unique connection $\nabla$ in
  the projective class, which preserves the conformally symplectic
  structure $\ell$. Then $\nabla$ defines a principal connection
  $\ga\in\Om^1(\Cal G_0,\frak g_0)$. Now we can pull back this
  connection form to $\Cal G_0^\#$ and then restrict it to $T^{-1}\Cal
  G_0^\#$. The result is a $G_0$--equivariant section of $L(T^{-1}\Cal
  G_0^\#,\frak g_0)$ which by construction reproduces the generators
  of fundamental vector fields. Analogously to the standard
  construction of induced connections, this gives rise to a partial
  contact connection $\nabla^\#$ on $H\subset TM^\#$.

  By construction, the local flows $\Fl^{\xi_0}_t$ satisfy
  $q_0\o\Fl^{\xi_0}_t=q_0$ and they preserve the contact distribution
  $H\subset TM^\#$. This immediately implies that these flows
  preserve the pullback of $\ga$ and its restriction to the contact
  distribution and hence the partial contact connection
  $\nabla^\#$. Now let $\al$ be the unique contact form on $M^\#$ such
  that $\al(\xi)=1$ and $\om$ the unique symplectic form on $M$ in the
  given conformal symplectic class such that $q^*\om=d\al$.  Consider
  sections $\eta_1$, $\eta_2$ and $\zeta$ of $H\subset TM^\#$ which
  project onto vector fields $\underline{\eta}_1$,
  $\underline{\eta}_2$ and $\underline{\zeta}$ on $M$. Then by
  definition, the Levi bracket $\Cal L(\eta_1,\eta_2)$ is equal to
  $\al([\eta_1,\eta_2])(\xi+H)=-(\om(\underline{\eta}_1,\underline{\eta}_2)\o
  q)(\xi+H)$. Now by definition of the partial connection on $Q$
  induced by $\nabla^\#$, we get
$$
\nabla^Q_\zeta\Cal L(\eta_1,\eta_2)=\Cal L(\nabla^\#_\zeta
\eta_1,\eta_2)+\Cal L(\eta_1,\nabla^\#_\zeta \eta_2). 
$$ 
By construction $\nabla^\#_\zeta \eta_i$ is a lift of
$\nabla_{\underline{\zeta}}\underline{\eta}_i$ for $i=1,2$. From
Proposition \ref{prop3.1}, we then know that $\nabla\omega=0$ and
inserting this, we easily conclude that $\nabla^Q(\xi+H)=0$. Since
$\xi$ is an infinitesimal contactomorphism, this implies that the
projection $\pi:TM^\#\to H$ induced by $\nabla^\#$ satisfies
$\pi(\xi)=0$, which of course determines $\pi$ completely. 

Using these facts, we can now continue completely parallel to the
proof of Theorem \ref{thm2.9}. We see that
$\pi([\eta_1,\eta_2])\in\Ga(H)$ is the unique section of $H$ lifting
$[\underline{\eta}_1,\underline{\eta}_2]$. Using this,
torsion--freeness of $\nabla$ implies that $\nabla^\#$ has vanishing
contact torsion. Hence we can use its contact projective equivalence
class to define a contact projective structure on $M$, and since the
flows of $\xi$ even preserve the partial contact connection
$\nabla^\#$, it is an infinitesimal automorphism of this structure.

  (3) Let $\ph:M\to\tilde M$ be the morphism of conformally Fedosov
structures and let $\ph^\#:M^\#\to\tilde M^\#$ be a contactomorphism
lifting $\ph$. Then the first step of the proof of Theorem \ref{thm2.8}
shows that $(\ph^\#)^*\tilde\xi=\la\xi$ for a nowhere vanishing,
locally constant function $\la$. Now there are natural lifts
$\Ph:\Cal G_0\to\tcg_0$ and $\Ph^\#:\Cal G_0^\#\to\tcg_0^\#$ to the
frame bundles and $(\Ph^\#)^*\tilde\xi_0=\la\xi_0$. In view of the
constructions of $\Ph^\#$ and $q_0$ and $\tilde q_0$, we see that
$\tilde q\o\ph^\#=\ph\o q$ implies $\tilde q_0\o\Ph^\#=\Ph\o q_0$.

  Now a morphism of conformally Fedosov structures preserves both the
  conformally symplectic structure and the projective class. Hence the
  uniqueness statement in statement (iii) in Proposition \ref{prop3.1}
  implies that it is also compatible with the (unique) connections in
  the projective class which preserve the conformally symplectic
  structures. Hence for their connection forms, we get
  $\Ph^*\tilde\ga=\ga$, which together with the above implies that
  $(\Ph^\#)^*\tilde q_0^*\tilde\ga=q_0^*\Ph^*\tilde\ga=q_0^*\ga$. Together
  with the fact that $(\Ph^\#)^*\tilde\xi_0=\la\xi_0$, this implies
  that $\Ph^\#$ is compatible with the partial contact connections
  constructed in (2) and hence in particular with the contact
  projective structures they generate. 
\end{proof}

Parallel to Theorem \ref{thm2.10}, there also is a relation between
conformally Fedosov structure and a class of special symplectic
connections. The relevant connections here are called
\textit{symplectic connections of Ricci type}, see \cite{BC} and
\cite{Cahen-Schwachhoefer}.  The module of formal curvatures of
symplectic connections splits over $\frak{sp}(2n)$ as a direct sum of
the kernel of the Ricci--type contractions and a complementary
submodule. A symplectic connection of Ricci type is a torsion--free
symplectic connection for which the curvature is concentrated in that
complementary submodule.

For a symplectic manifold $(M,\om)$ and a torsion--free symplectic
connection $\nabla$, the underlying conformally symplectic and
projective structures evidently define a conformally Fedosov
structure. On the other hand, for a general conformally Fedosov
structure, the connection $\nabla$ from part (iii) of Proposition
\ref{prop3.1} has curvature in $\La^2T^*M\otimes\frak{sp}(TM)$ by the last
part of the proposition. Hence the concept of being of Ricci type
makes sense for this connection.

\begin{cor}\label{cor3.3}
  The distinguished connection $\nabla$ of a conformally Fedosov
  structure from part (iii) of Proposition \ref{prop3.1} is of Ricci type
  if and only if any local parabolic contactification of the structure
  is locally flat as a contact projective structure.
\end{cor}
\begin{proof}
  Consider a local parabolic contactification $q:M^\#\to M$. In the
  construction of the contact projective structure on $M^\#$ in the
  proof of part (2) of Theorem \ref{thm3.3}, we used $\nabla$ to
  construct a partial contact connection $\nabla^\#$ on $M^\#$. The
  verifications in that proof actually show that $\nabla^\#$ is the
  partial contact connection coming from the Weyl structure determined
  by the infinitesimal automorphism $\xi$ giving rise to $q$.

The parabolic geometry determined by a contact projective structure
with vanishing contact torsion is well known to be torsion--free, see
Section 4.2.6 of \cite{book}. There it is also verified that the only
harmonic curvature of a contact projective structure has homogeneity 2
and is a section of the bundle induces by the highest weight subspace
in $\La^2\frak g_{-1}^*\otimes\frak g_0$, which is exactly the kernel
of the Ricci--type contraction. Now using Theorem 4.2.2 of
\cite{book}, one concludes as in the proof of Theorem \ref{thm2.10} that 
 this harmonic curvature can be computed from $\nabla^\#$. Using the
 information from the proof of part (2) of Theorem \ref{thm3.3} one
 deduces as in the proof of Theorem \ref{thm2.10} that vanishing of the
 component is equivalent to vanishing of the corresponding component
 of the curvature of $\nabla$, which implies the result. 
\end{proof}

\subsection{Example}\label{3.4} 
There is an analog of the global contactifications discussed in
Section \ref{2.6} in the context of conformally Fedosov
structures. This actually is the example relevant for the results on
integral geometry in \cite{E-G} which motivate the developments in this
series of articles. 

Let us realize the standard symplectic form in dimension $2n+2$ as the
imaginary part of the standard Hermitian form on $\Bbb C^{n+1}$. This
gives rise to an action of $G=Sp(2n+2,\Bbb R)$ on $\Bbb C^{n+1}$
(which does not preserve the complex structure). Viewing $S^{2n+1}$ as
the space of real rays in $\Bbb C^{n+1}$ we get a transitive action of
$G$, which makes $S^{2n+1}$ into the homogeneous model for contact
projective structures, compare with Section 4.2.6 of \cite{book}. Now
the standard action of $U(1)$ on $\Bbb C^{n+1}$ by complex
multiplication defines a subgroup in $G$ and hence an action of $U(1)$
on $S^{2n+1}$ by automorphisms of the contact projective
structure. The tangent space of $S^{2n+1}$ in a point corresponding to
a ray in $\Bbb C^{n+1}$ can be viewed as the quotient of the real
vector space $\Bbb C^{n+1}$ by the line spanned by the ray. The
contact subspaces then correspond to the quotient of the symplectic
orthocomplement of the ray by this line. This readily implies that the
infinitesimal generator of our $U(1)$--action is transversal
everywhere. Of course the $U(1)$--action on $S^{2n+1}$ admits a global
quotient, namely the canonical projection $q:S^{2n+1}\to\Bbb CP^n$
mapping a real ray in $\Bbb C^{n+1}$ to the complex line it spans.

Hence we obtain exactly the same quotient map as in the first example
in Section \ref{2.6}, where we started from the CR--structure on
$S^{2n+1}$, which is also preserved by the $U(1)$--action. It is
actually easy to see that the Weyl connections associated to the
infinitesimal generator of the $U(1)$--action via the CR--structure
and via the projective contact structure coincide. But this implies
that the underlying conformally Fedosov structure on $\Bbb CP^n$
actually is the one determined by the Levi--Civita connection of its
K\"ahler metric.

\begin{remark}\label{rem3.4}
  The last bit of this example actually exhibits a general
  phenomenon. Suppose that we have a torsion--free PCS--structure (not
  of type $C_n$) on a smooth manifold $M$. Then via Proposition
  \ref{prop3.1}, the (torsion--free) canonical connection actually
  determines an underlying conformally Fedosov structure. This may
  look surprising, since there is no corresponding construction for an
  ``underlying contact projective structure'' associated to a
  (torsion--free) parabolic contact structure of different
  type. However, such an underlying structure does become available in
  the presence of a transversal infinitesimal automorphism. Indeed,
  one can then consider the Weyl connection determined by this
  transversal infinitesimal automorphism and restrict it to a partial
  contact connection. General results then ensure that in the case of
  a torsion--free parabolic contact structure, this partial contact
  connection has vanishing contact torsion and hence determines a
  torsion--free contact projective structure, which is preserved by the
  given infinitesimal automorphism. Hence a contactification of the
  initial PCS--structure at the same time defines a contactification
  of the underlying conformally Fedosov structure.
\end{remark}

\end{document}